\RequirePackage{fix-cm}
\documentclass[smallextended]{svjour3}       
\smartqed  
\usepackage{graphicx}
\usepackage{amsfonts}
\usepackage{amsmath}
\usepackage{amssymb}
\usepackage{graphicx}
\usepackage{caption}

\def\Frac#1#2{\frac{\displaystyle{#1}}{\displaystyle{#2}}}

\begin{document}

\title{Uniform asymptotic expansions for Laguerre polynomials and related confluent
hypergeometric functions 
\thanks{The authors acknowledge support from \emph{Ministerio de Econom\'{\i}a y
Competitividad}, project MTM2015-67142-P (MINECO/FEDER, UE).}}

\titlerunning{Uniform asymptotic expansions for Laguerre polynomials}

\author{T. M. Dunster
\and A. Gil
\and J. Segura }

\institute{T. M. Dunster \at
              Department of Mathematics and Statistics\\
              San Diego State University. 5500 Campanile Drive San Diego, CA, USA.\\
              \email{mdunster@mail.sdsu.edu}   
           \and
           A. Gil \at
              Departamento de Matem\'atica Aplicada y CC. de la Computaci\'on.\\
ETSI Caminos. Universidad de Cantabria. 39005-Santander, Spain.\\
              \email{amparo.gil@unican.es} 
\and
           J. Segura \at
              Departamento de Matem\'aticas, Estad\'{\i}stica y Computaci\'on.\\
            Universidad de Cantabria. 39005-Santander, Spain.\\
              \email{segurajj@unican.es}          
}

\date{Received: date / Accepted: date}
\maketitle

\begin{abstract}
Uniform asymptotic expansions involving exponential and Airy functions are
obtained for Laguerre polynomials $L_{n}^{(\alpha)}(x)$, as well as complementary 
confluent hypergeometric functions. The expansions are valid for $n$ large and 
$\alpha$ small or large, uniformly
for unbounded real and complex values of $x$. 
The new expansions extend the range of 
computability of $L_n^{(\alpha)}(x)$ compared to previous expansions, 
in particular with respect to higher terms and large values of $\alpha$. 
Numerical evidence of their accuracy for real and complex values of $x$ is provided.


\end{abstract}

\numberwithin{equation}{section}

\section{Introduction}

\label{intro}

In this paper we shall obtain computable asymptotic expansions for Laguerre
polynomials $L_{n}^{(\alpha)}(x)$, for the case $n$ large and ${\alpha}$ small
or large. These are defined by
\begin{equation}
L_{n}^{(\alpha)}(x)=\sum_{k=0}^{n}\left(
\begin{array}
[c]{c}%
n+\alpha\\
n-k
\end{array}
\right)  \frac{(-x)^{k}}{k!}.\label{eq0}%
\end{equation}
In terms of the confluent hypergeometric functions we have
\begin{equation}
L_{n}^{(\alpha)}(x)=\frac{\Gamma(n+\alpha+1)}{n!}\mathbf{M}(-n,\alpha+1,x)
=\frac{(-1)^{n}}{n!}U(-n,\alpha+1,x),\label{eqLMU}%
\end{equation}
where $\mathbf{M}$ denotes Olver's scaled confluent hypergeometric function
\cite[Chap. 7, sec. 9]{Olver:1997:ASF}.

We will use the techniques described in \cite{Dunster:2017:COA} for computing
uniform asymptotic expansions of turning point problems. These involve Airy
function expansions for solutions of differential equations having a simple
turning point, but in a form that differs from the classical Airy function
expansions of \cite[Chap. 11]{Olver:1997:ASF}. Specifically, the coefficients
in these new expansions are significantly easier to compute than previously,
and close to the turning point the method utilises Cauchy's integral formula.
We remark that the Cauchy integral method of \cite{Dunster:2017:COA} has
potential applications to other forms of differential equations, and also the
method was subsequently used in \cite{SAPM:SAPM12172} to compute coefficients
appearing in certain asymptotic expansions of integrals.

In this paper we shall also obtain Liouville-Green (L-G) expansions for
Laguerre polynomials and confluent hypergeometric functions in domains that do
not contain turning points, and these involve the exponential function. The
form of these expansions differs from the standard ones (see \cite[Chap.
10]{Olver:1997:ASF}), insomuch the coefficients in the expansions appear in
the argument of the exponential. The advantage of this form is that the
coefficients are again generally easier to compute. Indeed, L-G expansions of
this form were used in \cite{Dunster:2017:COA} to obtain the new form of the
Airy expansions described above.

A recent reference on the computation of Laguerre polynomials using
asymptotics is \cite{Gil:2017:ECL}, where three different asymptotic
approximations are used: two expansions in terms of Bessel functions (from
\cite{Frenzen:1988:UAE} and \cite[Sect. 10.3.4]{Temme:2015:AMF}) and the Airy
expansion of Frenzen and Wong \cite{Frenzen:1988:UAE}. These expansions, and
in particular the Airy expansion, give accurate approximations for large $n$
when $\alpha$ is not large. For large $\alpha$ some expansions are available,
like for instance the expansions for Whittaker functions of
\cite{Dunster:1989:UAE} and \cite{Olver:1980:WFW}, or some of the expansions
in \cite{Temme:1990:AEF}, but the coefficients of these expansions are very
hard to compute. In this paper, we provide expansions that are computable and
which can also be used for large $\alpha$, and we provide numerical evidence
of their accuracy for real and complex values of $x$.

The Laguerre polynomials satisfy the following form of the confluent
hypergeometric equation
\begin{equation}
x\frac{d^{2} y}{dx^{2}}+(\alpha+1-x)\frac{dy}{dx}+ny=0.\label{eq1}%
\end{equation}
By the transformation $\tilde{y}=x^{(\alpha+1)/2}\exp\left( -\frac{1}%
{2}x\right) y$ we can remove the first derivative in the usual manner, and we
get
\begin{equation}
\frac{d^{2}\tilde{y}}{dx^{2}}=\frac{x^{2}-2(2n+\alpha+1)x+\alpha^{2}-1}%
{4x^{2}}\tilde{y}.\label{eq2}%
\end{equation}
On replacing $x$ by $uz$, where $u=n+\frac{1}{2}$ (a choice being explained
after (\ref{eq13a}) below), this can then be re-written in the form
\begin{equation}
\frac{d^{2} w}{dz^{2}}=\left\{ u^{2} f(a,z)+g(z)\right\} w,\label{eq3}%
\end{equation}
with a solution $w=z^{(\alpha+1)/2}\exp\left( -\frac{1}{2} uz\right)
L_{n}^{(\alpha)}(uz)$. On defining $a$ by
\begin{equation}
\alpha=u\left( a^{2}-1\right)  \ (a\geq0),\label{eq4}%
\end{equation}
we have in (\ref{eq3})
\begin{equation}
f(z) =\frac{(z-z_{1})(z-z_{2})}{4z^{2}},\ g(z) =-\frac{1}{4z^{2}},\label{eq5}%
\end{equation}
where
\begin{equation}
z_{1}=(a-1)^{2},\ z_{2}=(a+1)^{2}.\label{eq6}%
\end{equation}
The latter are the turning points for large $u$, which is assumed here. These
turning points coalesce with each other when $a=0$ ($\alpha=-u)$, and $z_{1}$
coalesces with the pole at $z=0$ when $a=1$ ($\alpha=0$).

We shall consider the following cases separately.

Case 1a: We obtain expansions for $z$ lying in domains containing the turning
point $z_{1}$, which includes the interval $0\leq z\leq z_{2} -\delta$ (here
and throughout $\delta>0$). The turning point $z=z_{2}$ is excluded, and this
is covered in Case 2 below. We assume $1+\delta\leq a^{2}\leq a_{1}<\infty$
for fixed $a_{1}\in(1,\infty)$. Thus $z_{1}$ cannot coalesce with the pole at
$z=0$ nor with the other turning point $z=z_{2}$, and $\alpha$ is positive and
large, satisfying $u\delta\leq\alpha\leq u(a_{1} -1)$.

Case 1b: This is the same as Case 1a, except now we consider $\alpha$
negative. In particular, we assume $0<a_{0}\leq a^{2}\leq1-\delta$ (for fixed
$a_{0}\in(0,1)$), and in this case $-u(1-a_{0}) \leq\alpha\leq-u\delta<0$.
Again $z_{1}$ cannot coalesce with the pole at $z=0$ nor with the other
turning point $z=z_{2}$.

Case 2: In this case expansions are derived for $z$ lying in domains
containing the turning point $z_{2}$, including the interval $z_{1}+\delta\leq
z<\infty$, but not $z_{1}$ or $0$. Here $a$ lies in a larger interval than
Cases 1a and 1b, namely $0<a_{0}\leq a^{2}\leq a_{1}<\infty$, and hence
$-u(1-a_{0})\leq\alpha\leq u(a_{1}-1)$. Note this can include large negative
values of $\alpha$ if $a_{0}<1$. The turning points again cannot coalesce
($a=0$) in this case, but $z_{1}$ can coalesce with $z=0$ (since both points
are excluded).

The more general case $0\leq a<\infty$ with $0\leq z<\infty$ will require an
application of asymptotic expansions valid for a coalescing turning point and
double pole, and also for two coalescing turning points. This will be studied
in a subsequent paper.

The plan of this paper is as follows. In section 2 we consider Case 1a and
construct L-G expansions for the functions. In this section a detailed
description of the Liouvile transformation is provided. In section 3 we obtain
Airy function expansions for Case 1a. In section 4 we obtain L-G and Airy
expansions for Case 1b, and in section 5 we likewise do this for Case 2.
Finally, in section 6 we present numerical results for the expansions of Cases
1a and 2 for Laguerre polynomials.

\section{L-G expansions: Case 1a}

We make the Liouville transformation
\begin{equation}
\xi=\int_{z}^{z_{1}}f^{1/2}(t) dt. \label{eq7}%
\end{equation}
On explicit integration (and noting that $\sqrt{z_{1} z_{2}}=\left\vert
a^{2}-1\right\vert $) we obtain
\begin{equation}%
\begin{array}
[c]{ll}%
\xi= & \dfrac{1}{2}\left(  a^{2}+1\right)  \ln\left\{  a^{2}+1-z-S\left(
z\right)  \right\}  -\dfrac{1}{2}S(z) -\max\left\{  a^{2},1\right\}  \ln(2a)\\
& \\
& +\dfrac{1}{2}\left\vert a^{2}-1\right\vert \ln\left\{  \dfrac{\left(
a^{2}-1\right)  ^{2}+\left\vert a^{2}-1 \right\vert S(z) -\left(
a^{2}+1\right)  z}{z}\right\}  .
\end{array}
\label{eq8}%
\end{equation}
In this $S(z)$ is given by
\begin{equation}
S(z) =\left\{  (z_{1}-z)(z_{2}-z)\right\}  ^{1/2}= \left\{  z^{2}-2\left(
a^{2}+1\right)  z+\left(  a^{2}-1\right)  ^{2}\right\}  ^{1/2}, \label{eq9}%
\end{equation}
where the principal logarithms are taken, and the branches of the square roots
are chosen so that $S(z)>0$ for $-\infty<z<z_{1}$, and is continuous elsewhere
in the plane having a cut on the interval $[z_{1},z_{2}] $. This means $S(z)
<0$ for $z_{2}<z<\infty$.

Thus $\xi$ is real and positive when $z\in(0,z_{1})$, and varies continuously
in the $z$ plane having cuts along $(-\infty,0]$ and $[z_{1},\infty) $. Note
that $z=z_{1}$ is mapped to $\xi=0$, and $z=z_{2}\pm i0$ (respectively above
and below the cut) is mapped to $\xi=\pm\min\left\{  a^{2},1\right\}  \pi i$.

\setkeys{Gin}{width=0.8\textwidth,height=0.8\textheight,keepaspectratio}

\begin{figure}[ptb]
\centering
\captionsetup{justification=centering} \vspace*{0.8cm}
\centerline{\includegraphics[height=6cm,width=8cm]{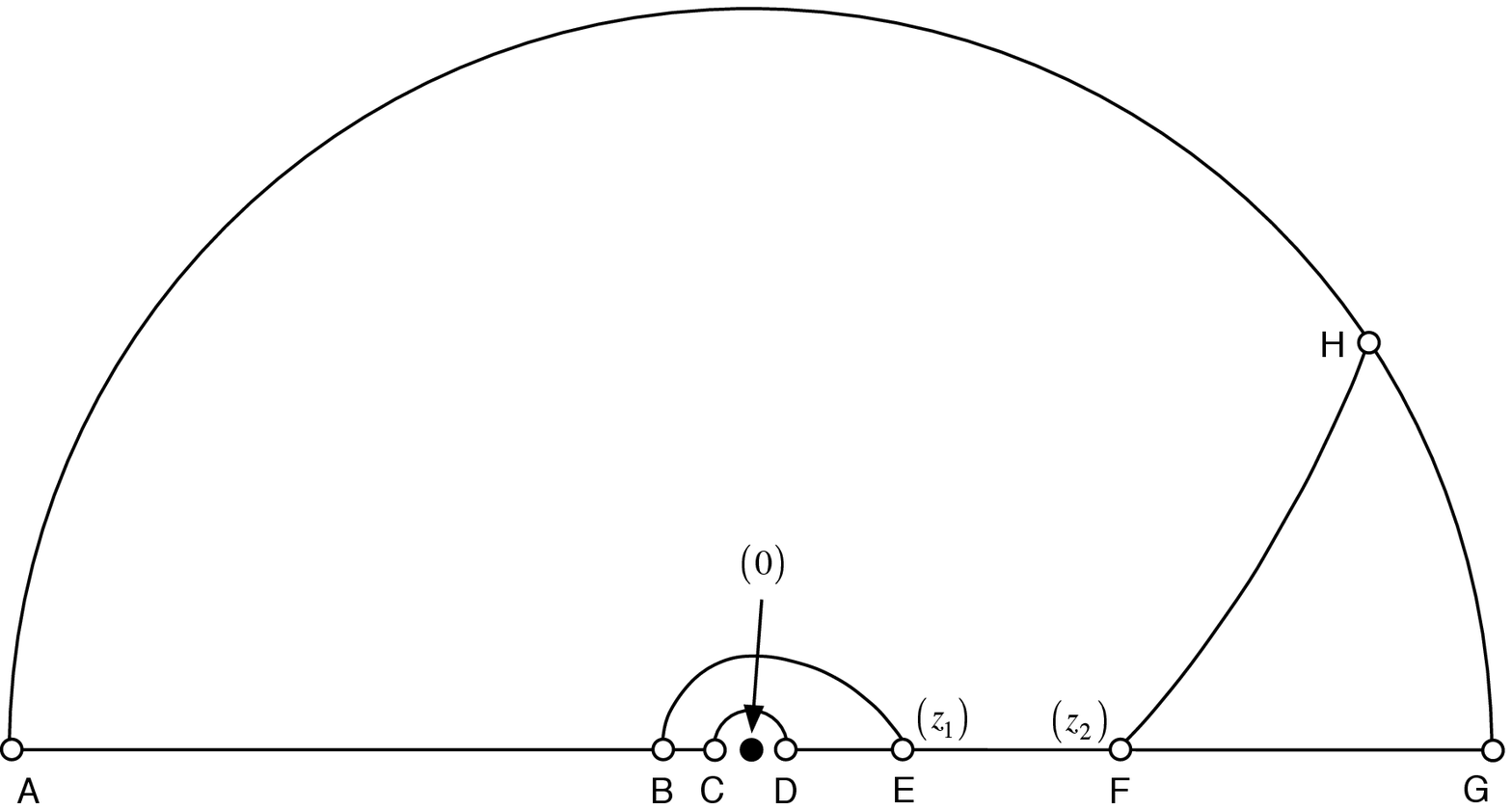}}\caption{$z$
plane.}%
\label{fig1zxi}%
\end{figure}

\begin{figure}[ptb]
\vspace*{0.8cm}
\centerline{\includegraphics[height=6cm,width=8cm]{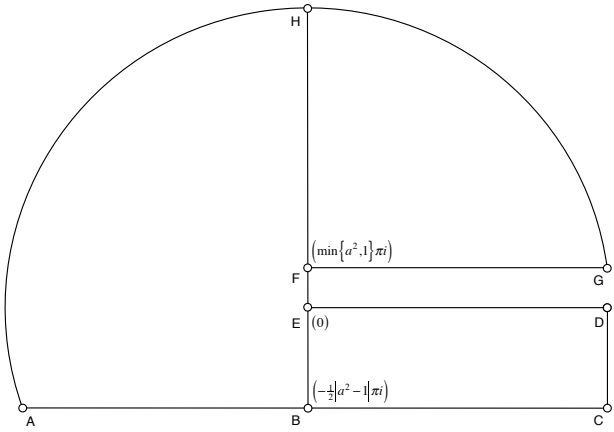}}\caption{$\xi$
plane.}%
\label{fig2zxi}%
\end{figure}

We find that $\xi\rightarrow+\infty$ as $z\rightarrow0^{+}$, such that
\begin{equation}
\xi=\tfrac{1}{2}\left\vert a^{2}-1\right\vert \left\{  2\ln\left\vert
{a^{2}-1}\right\vert -\ln\left(  z\right)  -1\right\}  -a^{2}\left\vert
\ln(a)\right\vert +\mathcal{O}(z). \label{eq10}%
\end{equation}

Furthermore, we have that
\begin{equation}
\xi=\tfrac{1}{2}z-\tfrac{1}{2}\left(  a^{2}+1\right)  \left\{  \ln
(z)+1\right\}  +a^{2}\ln(a) +\min\left\{  a^{2},1\right\}  \pi i
+\mathcal{O}\left(  {z^{-1}}\right)  , \label{eq51b}%
\end{equation}
as $z\rightarrow\infty$ in upper half plane.

Figs. \ref{fig1zxi} and \ref{fig2zxi} depict the $z$-$\xi$ map for $0\leq
\arg(z) \leq\pi$, with corresponding points labeled with the same letters. For
the lower half $z$ plane we can use from the Schwarz reflection principle that
$\xi(z) =\overline{\xi\left(  \bar{z}\right)  }$.

With (\ref{eq7}) and the new dependent variable given by
\begin{equation}
V=f^{1/4}(z) w, \label{eq11}%
\end{equation}
the differential equation (\ref{eq3}) is transformed to
\begin{equation}
\Frac{d^2 V}{d\xi^2}=\left\{  u^{2}+\phi(\xi)\right\}  V, \label{eq12}%
\end{equation}
where
\begin{equation}%
\begin{array}
[c]{ll}%
\phi\left(  \xi\right)  & = \dfrac{4f\left(  z\right)  {f}^{\prime\prime
}\left(  z\right)  -5{f}^{\prime2}\left(  z\right)  }{16f^{3}\left(  z\right)
}+\dfrac{g\left(  z\right)  }{f\left(  z\right)  }\\
& =-\dfrac{z\left\{  {4z^{3}-4\left(  {3a^{2}-1}\right)  \left(  {a^{2}%
-3}\right)  z+8\left(  {a^{2}+1}\right)  \left(  {a^{2}-1}\right)  ^{2}%
}\right\}  }{4\left(  {z-z_{1}}\right)  ^{3}\left(  {z-z_{2}}\right)  ^{3}}.
\end{array}
\label{eq13}%
\end{equation}

The function $\phi(\xi)$ is analytic in an unbounded domain $\Delta$ (say)
which excludes $\xi=0$ and $\xi=\pm\min\left\{  a^{2},1\right\}  \pi i$ (the
singularities corresponding to the turning points). Then the part of $\Delta$
corresponding to $0\leq\arg(z)\leq\pi$ is the entire region depicted in Fig.
\ref{fig2zxi}, except with the points $\xi=0$ and $\xi=\min\left\{
a^{2},1\right\}  \pi i$ (labeled E and F, respectively) excluded.

An important property is that
\begin{equation}
\phi(\xi) =\mathcal{O}\left(  \xi^{-2}\right)  , \label{eq13a}%
\end{equation}
as $\xi\rightarrow\infty$ in $\Delta$. We remark that the choice of the large
parameter in the form $u=n+\tfrac{1}{2}$, and the subsequent partitioning of
(\ref{eq3}), resulted in the desired behaviour (\ref{eq13a}).

Asymptotic solutions of (\ref{eq12}), accompanied by explicit error bounds,
are given in \cite{Dunster:1998:AOT} by
\begin{equation}
V_{n,1}(u,\xi) =\exp\left\{  u\xi+\sum\limits_{s=1}^{n-1}\frac{E_{s}(\xi
)}{u^{s}}\right\}  +\varepsilon_{n,1}(u,\xi), \label{eq14}%
\end{equation}
and
\begin{equation}
V_{n,2}(u,\xi) =\exp\left\{  -u\xi+\sum\limits_{s=1}^{n-1}(-1)^{s}\frac
{E_{s}(\xi)}{u^{s}}\right\}  +\varepsilon_{n,2}(u,\xi), \label{eq14a}%
\end{equation}
where the coefficients $E_{s}(\xi)$ will be defined below.

In (\ref{eq14}) $\varepsilon_{n,1}(u,\xi)=e^{u\xi} \mathcal{O}\left(
u^{-n}\right)  $ uniformly in a certain domain $\Xi_{1}$, with $e^{-u\xi
}\varepsilon_{n,1}(u,\xi)\rightarrow0$ as $\mathrm{Re} \xi\rightarrow-\infty$
in the domain, and likewise in (\ref{eq14a}) $\varepsilon_{n,2}(u,\xi)
=e^{-u\xi}\mathcal{O}(u^{-n})$ uniformly in a domain $\Xi_{2}$, with $e^{u\xi
}\varepsilon_{n,2}(u,\xi)\rightarrow0$ as $\mathrm{Re} \xi\rightarrow\infty$
in the domain. By virtue of (\ref{eq13a}), both $\Xi_{1}$ and $\Xi_{2}$ are
unbounded, and are defined as follows.

In general the domain of validity $\Xi_{1}$ comprises the $\xi$ point subset
of $\Delta$ for which there is a path $\mathcal{P}_{1}$ (say) linking $\xi$
with $\alpha_{1}=-\infty-\frac{1}{2}\left\vert a^{2}-1\right\vert \pi i$
(corresponding to $z=\infty e^{\pi i}$) and having the properties (i)
$\mathcal{P}_{1}$ consists of a finite chain of $R_{2}$ arcs (as defined in
\cite[Chap. 5, sec. 3.3]{Olver:1997:ASF}), and (ii) as $t$ passes along
$\mathcal{P}_{1}$ from $\alpha_{1}$ to $\xi$, $\mathrm{Re} (ut) $ is nondecreasing.

The domain $\Xi_{2}$ comprises the $\xi$ point subset of $\Delta$ for which
there is a path $\mathcal{P}_{2}$ (say) linking $\xi$ with $\alpha_{2}
=+\infty$ (corresponding to $z=0$) and having the properties (i)
$\mathcal{P}_{2}$ consists of a finite chain of $R_{2}$ arcs, and (ii) as $t$
passes along $\mathcal{P}_{2}$ from $\alpha_{2}$ to $\xi$, $\mathrm{Re}(ut)$
is nonincreasing.

Let $\Xi_{j}^{+}$ ($j=1,2$) denote the subsets of $\Xi_{j}$ corresponding to
$0\leq\arg(z) \leq\pi$. Then $\Xi_{1}^{+}$ coincides with the domain shown in
Fig. \ref{fig2zxi}, but with the points $\xi=0$ and $\xi=\min\left\{
a^{2},1\right\}  \pi i$ excluded. We denote by $D_{1}^{+}$ the $z$-domain
corresponding to $\Xi_{1}^{+}$, and this is the entire upper half plane
$\mathrm{Im}(z) \geq0$ excluding the turning points $z=z_{1},z_{2}$.

\begin{figure}[ptb]
\vspace*{0.8cm}
\centerline{\includegraphics[height=6cm,width=8cm]{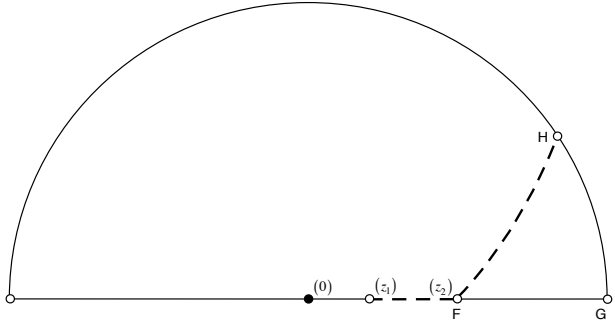}}\caption{Domain
$D_{2}^{+}$ in $z$ plane. The interval $\left[  z_{1},z_{2}\right]  $ and
boundary FH are excluded.}%
\label{fig3zxi}%
\end{figure}

On the other hand, due to the monotonicity requirement (ii), the domain
$\Xi_{2}^{+}$ must also exclude the unbounded region FGH of Fig.
\ref{fig2zxi}, as well as points on the segment EF. The corresponding
$z$-domain $D_{2}^{+}$ (say) is the unshaded region depicted in Fig.
\ref{fig3zxi}, where the interval $[z_{1},z_{2}]$ and the boundary curve FH
must be excluded; the curve FH extends from $z=z_{2}$ to $z=\infty$ and is
given parametrically by
\begin{equation}
\int_{z_{2}}^{z}{f^{1/2}\left(  t\right)  dt=\tau i,\ 0\leq\tau<\infty.}
\label{eq14b}%
\end{equation}
We note that $z=+\infty$ is not contained in $D_{2}^{+}$ (but is contained in
$D_{1}^{+}$).

Returning to the asymptotic expansions (\ref{eq14}) and (\ref{eq14a}), the
coefficients are given by
\begin{equation}
E_{s}(\xi) =\int F_{s}(\xi) d\xi\ (s=1,2,3,\cdots), \label{eq16}%
\end{equation}
where
\begin{equation}
F_{1} (\xi) =\frac{1}{2}\phi(\xi) ,\ F_{2}(\xi) =-\frac{1}{4}\phi(\xi),
\label{eq17}%
\end{equation}
and
\begin{equation}
F_{s+1}(\xi) =-\frac{1}{2}F_{s}^{\prime}(\xi) -\frac{1}{2}\sum\limits_{j=1}%
^{s-1} F_{j}(\xi) F_{s-j}(\xi)\,\, (s=2,3,\cdots). \label{eq18}%
\end{equation}
Primes are derivatives with respect to $\xi$.

The choice of integration constants in (\ref{eq16}) must be the same for both
of the L-G solutions (\ref{eq14}) and (\ref{eq14a}). As discussed in
\cite{Dunster:2017:COA}, the constants associated with the even terms ($s=2j$,
$j=1,2,3\cdots$) are arbitrary, but those for the odd terms ($s=2j+1$,
$j=0,1,2,\cdots$) must be precisely chosen, as described below.

It is preferable to work in terms of $z$. Using
\begin{equation}
\frac{d\xi}{dz}=f^{1/2}(z)=\frac{S(z)}{2z}, \label{eq19}%
\end{equation}
and writing $\hat{E}_{s}(z) =E_{s}\left(  \xi(z)\right)  $ and $\hat{F}%
_{s}(z)=F_{s}\left(  \xi(z)\right)  $, we have
\begin{equation}
\hat{E}_{s} (z) =\int\frac{\hat{F}_{s} (z) S (z)}{2z}dz \ (s=1,2,3,\cdots),
\label{eq20}%
\end{equation}
where
\begin{equation}
\hat{F}_{1} (z) =\frac{1}{2}\phi\left(  \xi(z)\right)  ,\ \hat{F}_{2} (z)
=-\frac{z}{2 S(z)} \frac{d\phi\left(  \xi(z)\right)  }{dz}, \label{eq21}%
\end{equation}
and
\begin{equation}
\hat{F}_{s+1}(z) =-\frac{z}{S(z)}\hat{F}_{s}^{\prime}(z) -\frac{1}{2}%
\sum\limits_{j=1}^{s-1} \hat{F}_{j} (z) \hat{F}_{s-j} (z)\ (s=2,3,\cdots).
\label{eq22}%
\end{equation}
We find by induction that
\begin{equation}
\frac{\hat{F}_{s} (z) S(z)}{2z}= \frac{R_{2s+1}(z)}{S^{3s+2}(z)}, \label{eq23}%
\end{equation}
where $R_{n}(z)$ is a polynomial of degree $n$. From this we can show that
\begin{equation}
\hat{E}_{2j}(z) =\frac{zT_{4j-1}^{(e)}(z)}{S^{6j}(z)}, \label{eq24}%
\end{equation}
and
\begin{equation}
\hat{E}_{2j+1} (z) =\frac{T_{6j+3}^{(o)}(z)}{S^{6j+3}(z)}, \label{eq25}%
\end{equation}
where $T_{n}^{(e)}(z)$ and $T_{n}^{(o)}(z)$ are also polynomials of degree
$n$, provided the integration constants in (\ref{eq20}) are chosen
appropriately (which we assume).

For example, for the coefficient $\hat{E}_{1}(z)$, we have
\begin{equation}%
\begin{array}
[c]{ll}%
T_{3}^{(o)}(z) = & \dfrac{1}{48a^{2}}\left[  \left(  a^{2}-1\right)
^{2}\left\{  \left(  a^{2}-1\right)  ^{2}-4a^{2}\right\}  -3\left(
a^{2}+1\right)  ^{3} z \right. \\
& \\
& \left.  +3\left(  a^{4}+6a^{2}+1\right)  z^{2} -\left(  a^{2}+1\right)
z^{3} \right]  .
\end{array}
\label{eqT3}%
\end{equation}

We note here that $\zeta^{1/2}\hat{E}_{2j+1}(z)$, regarded as a function of
$\zeta=\left(  \frac32\xi\right)  ^{2/3}$, is meromorphic at $\zeta=0$. From
\cite{Dunster:2017:COA} this is a requirement for the subsequent Airy function
expansions to be valid at the turning point.

We also note that
\begin{equation}
\hat{E}_{s}(\infty)=\lambda_{s}, \label{eq26}%
\end{equation}
(say), where $\lambda_{2j}=0$ ($j=1,2,3,\cdots$), whereas the odd terms are
non-zero. The first three of the odd terms are found to be
\begin{equation}
\lambda_{1} =-\frac{1+a^{2}}{48a^{2}},\lambda_{3} =\frac{7\left(
1+a^{6}\right)  } {5760 a^{6}},\lambda_{5}=-\frac{31\left(  1+a^{10}\right)
}{80640 a^{10}}. \label{eq27}%
\end{equation}

We similarly find
\begin{equation}
\hat{E}_{s}(0) =\mu_{s}, \label{eq28}%
\end{equation}
(say), where $\mu_{2j}=0$ ($j=1,2,3,\cdots$), and again the odd terms are
non-zero. The first two of these are found to be
\begin{equation}%
\begin{array}
[c]{l}%
\mu_{1}=\Frac{a^{4}-6a^{2}+1}{48a^{2}\left\vert a^{2}-1\right\vert },\\
\\
\mu_{3}%
=-\Frac{7a^{12}-21a^{10}+21a^{8}-30a^{6}+21a^{4}-21a^{2}+7} {5760a^{6}\left\vert a^{2}-1 \right\vert ^{3}}.
\end{array}
\label{eq29}%
\end{equation}

Recall in this case we are assuming that $1+\delta\leq a^{2}\leq a_{1}<\infty
$, for fixed $a_{1}\in( 1,\infty)$. Thus from (\ref{eq4}) $\alpha$ is
positive, and so using
\begin{equation}
L_{n}^{(\alpha)}(x)=\frac{\Gamma(n+\alpha+1)}{n!\Gamma(\alpha+1)} \left\{
1+\mathcal{O}(x)\right\}  \ (x\rightarrow0) , \label{eq30}%
\end{equation}
we therefore have the solution of (\ref{eq3}) given by
\begin{equation}
w_{0} (u,z) \equiv z^{(\alpha+1)/2} e^{-uz/2}L_{n}^{(\alpha)}(uz),
\label{eq31}%
\end{equation}
which has the unique behaviour as $z\rightarrow0$
\begin{equation}
w_{0} (u,z) \equiv\frac{\Gamma(n+\alpha+1)z^{(\alpha+1)/2}} {n!\Gamma
(\alpha+1)} \left\{  1+\mathcal{O}(z)\right\}  . \label{eq32}%
\end{equation}

We now match solutions that are recessive at $z=0$. From (\ref{eq11}) we have
by uniqueness (up to a multiplicative constant) of such solutions that%
\begin{equation}
w_{0} (u,z) \propto\dfrac{z^{1/2}V_{n,2}(u,\xi)} {\left\{  (z_{1}%
-z)(z_{2}-z)\right\}  ^{1/4}}. \label{eq32a}%
\end{equation}
Then, using (\ref{eq10}), (\ref{eq14a}), (\ref{eq28}), (\ref{eq31}) and
(\ref{eq32}), we arrive at
\begin{equation}%
\begin{array}
[c]{l}%
L_{n}^{(\alpha)}(uz) \sim\dfrac{\Gamma(n+\alpha+1)} {n!\Gamma(\alpha+1)}
\left(  \dfrac{\alpha}{u}\right)  ^{1/2}\left(  \dfrac{u}{u+\alpha}\right)
^{u/2} \left\{  \dfrac{\alpha^{2}}{(u+\alpha)uez}\right\}  ^{\alpha/2}\\
\times\dfrac{1}{\left\{  (z_{1}-z)(z_{2}-z)\right\}  ^{1/4}} \exp\left\{
\dfrac{1}{2}uz-u\xi+\displaystyle\sum_{s=1}^{\infty}(-1)^{s} \dfrac{\hat
{E}_{s} (z) -\mu_{s}}{u^{s}}\right\}  .
\end{array}
\label{eq33}%
\end{equation}
This is uniformly valid for $z\in D_{2}^{+}\cup D_{2}^{-}$, where $D_{2}^{-}$
is the conjugate region of $D_{2}^{+}$.

From \cite[Eq. 13.2.24]{Daalhuis:2010:CHF} a companion solution of (\ref{eq3})
is given by
\begin{equation}
w_{1}(u,z)\equiv(uz)^{(\alpha+1)/2}e^{uz/2}U(n+\alpha+1,\alpha+1,uze^{-\pi
i}).\label{eq50}%
\end{equation}
Now from \cite[Eq. 13.7.3]{Daalhuis:2010:CHF}
\begin{equation}
U(n+\alpha+1,\alpha+1,uze^{-\pi i})=\left(  uze^{-\pi i}\right)
^{-n-\alpha-1}\left\{  1+\mathcal{O}\left(  z^{-1}\right)  \right\}
\label{eq50b}%
\end{equation}
as $z\rightarrow\infty$ for $\left\vert \arg\left(  ze^{-\pi i}\right)
\right\vert \leq\frac{3}{2}\pi-\delta$, and we see that $w_{1}(u,z)$ has the
behaviour
\begin{equation}
w_{1}(u,z)=(-1)^{n+1}e^{-\alpha\pi i}(uz)^{-u-\alpha/2}e^{uz/2}\left\{
1+\mathcal{O}\left(  z^{-1}\right)  \right\}  ,\label{eq51}%
\end{equation}
and in particular this is the unique solution that is recessive in the
half-plane ${\left\vert {\arg\left(  {ze^{-\pi i}}\right)  }\right\vert
\leq{\frac{1}{2}}\pi}$.

Thus we similarly deduce by matching recessive solutions that
\begin{equation}
w_{1}(u,z)\propto\dfrac{z^{1/2}V_{n,1}(u,\xi)}{\left\{  (z_{1}-z)(z_{2}%
-z)\right\}  ^{1/4}},\label{eq51c}%
\end{equation}
and hence using (\ref{eq51b})%

\begin{equation}%
\begin{array}
[c]{l}%
U(n+\alpha+1,\alpha+1,uze^{-\pi i}) \sim-\dfrac{i\exp\left\{  u+\tfrac{1}%
{2}\alpha+\max(\alpha,0)\pi i\right\}  } {u^{(u+1)/2}(u+\alpha)^{(u+\alpha
)/2}(uz)^{\alpha/2}}\\
\\
\times\dfrac{1}{\left\{  (z_{1}-z)(z_{2}-z)\right\}  ^{1/4}} \exp\left\{
-\dfrac{1}{2}uz+u\xi+\displaystyle\sum_{s=1}^{\infty} \dfrac{\hat{E}_{s}(z)
-\lambda_{s}}{u^{s}}\right\}  .
\end{array}
\label{eq51a}%
\end{equation}

We remark that this expansion also holds for negative $\alpha$, and in
particular for $0<a_{0}\leq a^{2}\leq1-\delta$ (for fixed $a_{0}\in(0,1)$).
Likewise for its Airy expansion (\ref{eq43b}) below.

\section{Airy expansions: Case 1a}

We now obtain asymptotic expansions which are valid at $z=z_{1}$. These
involve Airy functions, and the standard Liouville transformation is given by
\begin{equation}
\zeta=\left(  \tfrac{3}{2}\xi\right)  ^{2/3},\ W=\zeta^{-1/4}f^{1/4}(z) w,
\label{eq34}%
\end{equation}
where $\xi$ is again given by (\ref{eq8}). Here $\zeta$ is defined to be
analytic at $z=z_{1}$ (see \cite[Chap. 11]{Olver:1997:ASF}), and moreover
$\zeta>0$ for $0<z<z_{1}$, and $\zeta<0$ for $z_{1}<z<z_{2}$.

With this transformation the differential equation (\ref{eq3}) takes the form
\begin{equation}
\dfrac{d^{2}W}{d\zeta^{2}}=\left\{  u^{2}\zeta+\psi(\zeta)\right\}  W,
\label{eq35}%
\end{equation}
where $\psi(\zeta)$ is analytic at $\zeta=0$ (i.e. $\xi=0$ and $z=z_{1}$), and
is given explicitly by
\begin{equation}
\psi(\zeta)=\tfrac{5}{16}\zeta^{-2}+\zeta\phi(\xi), \label{eq36}%
\end{equation}
where $\phi(\xi) $ is given by (\ref{eq13}).

From \cite{Dunster:2017:COA} the following three \textit{exact} solutions of
(\ref{eq35}) are given
\begin{equation}
W_{j}(u,\zeta) =\mathrm{Ai}_{j}\left(  u^{2/3}\zeta\right)  A(u,z)
+\mathrm{Ai}_{j}^{\prime}\left(  u^{2/3}\zeta\right)  B(u,z) \ (j=0,\pm1).
\label{eq37}%
\end{equation}
Here $\mathrm{Ai}_{j}(u^{2/3}\zeta)=\mathrm{Ai}(u^{2/3}\zeta e^{-2\pi ij/3})$,
which are the Airy functions that are recessive in the sectors
$\mathrm{\mathbf{S}}_{j}:=\left\{  \zeta:|\arg(\zeta e^{-2\pi ij/3})|\leq
\pi/3\right\}  $; see \cite[\S 9.2(iii)]{Olver:2010:AAR}. In our case, as
functions of $z$, $W_{0}(u,\zeta) $ is recessive at $z=0$, $W_{1}(u,\zeta)$ is
recessive at $z=\infty e^{\pi i}$, and $W_{-1}(u,\zeta) $ is recessive at
$z=\infty e^{-\pi i}$.

The coefficient functions $A(u,z)$ and $B(u,z)$ are analytic in a domain
containing $z=z_{1}$, and in \cite[Theorem 2.1]{Dunster:2017:COA} it was shown
that they possess the asymptotic expansions
\begin{equation}%
\begin{array}
[c]{l}%
A(u,z)\sim\exp\left\{  \displaystyle\sum_{s=1}^{\infty} \dfrac{\hat{E}%
_{2s}(z)+\tilde{a}_{2s}\xi^{-2s}/(2s)}{u^{2s}}\right\} \\
\times\cosh\left\{  \displaystyle\sum_{s=0}^{\infty} \dfrac{\hat{E}%
_{2s+1}(z)-\tilde{a}_{2s+1}\xi^{-2s-1}/(2s+1)}{u^{2s+1}}\right\}  ,
\end{array}
\label{eq38}%
\end{equation}
and
\begin{equation}%
\begin{array}
[c]{l}%
B(u,z)\sim\dfrac{1}{u^{1/3}\zeta^{1/2}}\exp\left\{  \displaystyle\sum
_{s=1}^{\infty} \dfrac{\hat{E}_{2s}(z)+a_{2s}\xi^{-2s}/(2s)}{u^{2s}}\right\}
\\
\times\sinh\left\{  \displaystyle\sum_{s=0}^{\infty} \dfrac{\hat{E}%
_{2s+1}(z)-a_{2s+1}\xi^{-2s-1}/(2s+1)}{u^{2s+1}}\right\}  ,
\end{array}
\label{eq39}%
\end{equation}
where $a_{1}=a_{2}=\frac{5}{72}$, $\tilde{a}_{1}=\tilde{a}_{2}=-\frac{7}{72}$,
and for $b=a$ and $b=\tilde{a}$
\begin{equation}
b_{s+1}=\frac{1}{2}(s+1)b_{s}+\frac{1}{2}\sum_{j=1}^{s-1}b_{j}b_{s-j}
\quad(s=2,3,\cdots). \label{eq39a}%
\end{equation}

\begin{figure}[ptb]
\vspace*{0.8cm}
\centerline{\includegraphics[height=6cm,width=8cm]{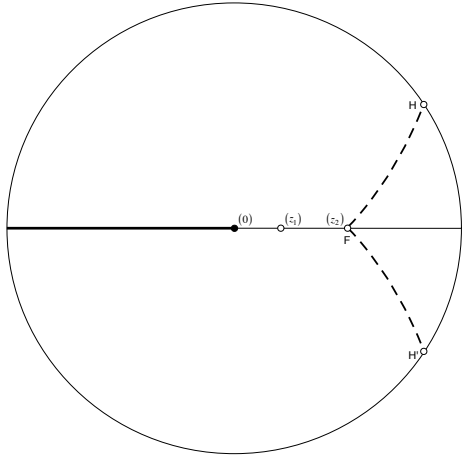}}\caption{Domain
$\mathbf{D}$ in $z$ plane.}%
\label{fig4zxi}%
\end{figure}

For our differential equation (\ref{eq35}) these expansions are uniformly
valid in an unbounded $z$ domain $\mathbf{D}$, which consists of all points
which can be linked to each of the singularities $z=0,z=\infty e^{\pm\pi i}$
by a "progressive" path: that is, a finite chain of $R_{2}$ arcs that does not
pass through $z_{2}$, and with the property that $\mathrm{Re}\xi(z)$ varies
monotonically as the path is traversed from one end to the other. From Fig.
\ref{fig2zxi} it is straightforward to show that $\mathbf{D}$ consists of the
intersection of $D_{1}^{+}$ and $D_{2}^{+}$ (i.e. $D_{2}^{+}$ itself), along
with the conjugate of $D_{2}^{+}$, and in addition all points on the interval
$[z_{1},z_{2})$.

Thus $\mathbf{D}$ is the unshaded region shown in Fig. \ref{fig4zxi}. In this
$\left\vert \arg(z)\right\vert \leq\pi$, the boundary curve FH emanating from
$z=z_{2}$ is defined by (\ref{eq14b}), and FH' is the conjugate curve. Points
on both curves (including $z=z_{2}$) are excluded from $\mathbf{D}$, but the
singularities $z=0$ and $z=\infty e^{\pm\pi i}$ lie in $\mathbf{D}$, as well
as of course the turning point $z=z_{1}$.

For each of the three solutions the respective domain of validity extends
beyond $\mathbf{D}$. For example, for $W_{0}(u,\zeta)$ certain points
accessible by crossing above or below the cut on the negative real axis can be
included. For $W_{1}(u,\zeta)$ certain points crossing above this cut can be
included, as well as points crossing the curve FH (but not FH'). For our
purposes the common domain of validity $\mathbf{D}$ suffices, since an
extension to a domain containing $z=z_{2}$ and $z=+\infty$ will be considered
in the next section, and for the confluent hypergeometric function analytic
continuation formulas can be used for values of $\arg(z)$ outside $[-\pi,\pi]$.

We now match the Airy function expansions with the Laguerre polynomial and
confluent hypergeometric function. Firstly we have, on matching the solutions
recessive at $z=0$, namely $W_{0}(u,\zeta)$ and $\zeta^{-1/4}f^{1/4}(z) w_{0}
(u,z)$ (see (\ref{eq31}) and (\ref{eq37})), that
\begin{equation}
W_{0}(u,\zeta) =C_{0}(u) \left\{  \dfrac{(z_{1}-z)(z_{2}-z)}{\zeta}\right\}
^{1/4} z^{\alpha/2}e^{-uz/2}L_{n}^{(\alpha)}(uz), \label{eq40}%
\end{equation}
for some constant $C_{0}(u)$.

Now as $z\rightarrow0^{+}$ ($\zeta\rightarrow+\infty$) we find from
(\ref{eq34}), (\ref{eq38}), (\ref{eq39}), and the behaviour of the Airy
function and its derivative for large argument (see \cite[\S 9.7(ii)]%
{Olver:2010:AAR}), that the LHS behaves as
\begin{equation}
W_{0}(u,\zeta)\sim\dfrac{1}{2\pi^{1/2}u^{1/6}\zeta^{1/4}}\exp\left\{
-u\xi-\sum\limits_{s=0}^{\infty}{\dfrac{\mu_{2s+1}}{u^{2s+1}}}\right\}  .
\label{eq41}%
\end{equation}
On the other hand, from (\ref{eq30}), the RHS of (\ref{eq40}) has the
behaviour
\begin{equation}%
\begin{array}
[c]{l}%
C_{0}\left(  u\right)  \left\{  {\dfrac{\left(  {z-z_{1}}\right)  \left(
{z-z_{2}}\right)  }{\zeta}}\right\}  ^{1/4}z^{\alpha/2}e^{-uz/2}L_{n}^{\left(
\alpha\right)  }\left(  {uz}\right) \\
\sim\dfrac{C_{0}\left(  u\right)  \Gamma\left(  {n+\alpha+1}\right)  \left(
{a^{2}-1}\right)  ^{1/2}z^{\alpha/2}}{n!\Gamma\left(  {\alpha+1}\right)
\zeta^{1/4}}.
\end{array}
\label{eq42}%
\end{equation}
Solving for $C_{0}\left(  u\right)  $ from (\ref{eq40}) - (\ref{eq42}), and on
referring to (\ref{eq4}) and (\ref{eq51b}), we then arrive at
\begin{equation}%
\begin{array}
[c]{l}%
L_{n}^{\left(  \alpha\right)  }\left(  {uz}\right)  \sim\dfrac{2\pi
^{1/2}u^{1/6}\Gamma\left(  {n+\alpha+1}\right)  }{n!\Gamma\left(  {\alpha
+1}\right)  }\left(  {\dfrac{\alpha}{u}}\right)  ^{1/2}\left(  {\dfrac
{u}{u+\alpha}}\right)  ^{u/2}\\
\\
\times\left(  {\dfrac{\alpha^{2}}{\left(  {u+\alpha}\right)  uez}}\right)
^{\alpha/2}\left\{  {\dfrac{\zeta}{\left(  {z_{1}-z}\right)  \left(  {z_{2}%
-z}\right)  }}\right\}  ^{1/4}\exp\left\{  {\dfrac{1}{2}uz+\sum\limits_{s=0}%
^{\infty}{\dfrac{\mu_{2s+1}}{u^{2s+1}}}}\right\} \\
\\
\times\left\{  {\mathrm{Ai}\left(  {u^{2/3}\zeta}\right)  A\left(
{u,z}\right)  +}\mathrm{Ai}^{\prime}{\left(  {u^{2/3}\zeta}\right)  B\left(
{u,z}\right)  }\right\}  .
\end{array}
\label{eq43}%
\end{equation}

Similarly, on matching solutions recessive at $z=\infty e^{\pi i}$ we assert
the existence of a constant $C_{1}\left(  u\right)  $ such that
\begin{equation}%
\begin{array}
[c]{l}%
W_{1}\left(  {u,\zeta}\right)  =C_{1}\left(  u\right)  \left\{  {\dfrac
{\left(  {z_{1}-z}\right)  \left(  {z_{2}-z}\right)  }{\zeta}}\right\}
^{1/4}\\
\times z^{\alpha/2}e^{uz/2}U\left(  {n+\alpha+1,\alpha+1,uze^{-\pi i}}\right)
.
\end{array}
\label{eq43aa}%
\end{equation}
The constant is found by comparing both sides as $z\rightarrow\infty e^{\pi
i}$, and as a result, on using
\begin{equation}%
\begin{array}
[c]{l}%
\mathrm{Ai}_{1}\left(  {u^{2/3}\zeta}\right)  A\left(  {u,z}\right)
+\mathrm{Ai}_{1}^{\prime}\left(  {u^{2/3}\zeta}\right)  B\left(  {u,z}\right)
\\
\sim\dfrac{e^{\pi i/6}}{2\pi^{1/2}u^{1/6}\zeta^{1/4}}\exp\left\{  {u\xi
+\sum\limits_{s=0}^{\infty}{\dfrac{\lambda_{2s+1}}{u^{2s+1}}}}\right\}  ,
\end{array}
\label{eq43a}%
\end{equation}
we find that%
\begin{equation}%
\begin{array}
[c]{l}%
U\left(  {n+\alpha+1,\alpha+1,uze^{-\pi i}}\right) \\
\\
\sim-\dfrac{2\pi^{1/2}e^{\pi i/3}}{u^{1/3}\left\{  {u\left(  {u+\alpha
}\right)  }\right\}  ^{\left(  {u+\alpha}\right)  /2}z^{\alpha/2}}\left\{
{\dfrac{\zeta}{\left(  {z_{1}-z}\right)  \left(  {z_{2}-z}\right)  }}\right\}
^{1/4}\\
\\
\times\exp\left\{  {-\dfrac{1}{2}uz+u+\dfrac{1}{2}\alpha+\max}\left(
{\alpha,0}\right)  {\pi i-\sum\limits_{s=0}^{\infty}{\dfrac{\lambda_{2s+1}%
}{u^{2s+1}}}}\right\} \\
\\
\times\left\{  {\mathrm{Ai}_{1}\left(  {u^{2/3}\zeta}\right)  A\left(
{u,z}\right)  +}\mathrm{Ai}_{1}^{\prime}{\left(  {u^{2/3}\zeta}\right)
B\left(  {u,z}\right)  }\right\}  .
\end{array}
\label{eq43b}%
\end{equation}

By the matching of solutions that are recessive as $z\rightarrow\infty e^{-\pi
i}$, we obtain corresponding L-G and Airy asymptotic expansions which are
given by (\ref{eq51a}) and (\ref{eq43b}) respectively, with $i$ replaced by
$-i$ in both, and ${\mathrm{Ai}_{1}}$ replaced by ${\mathrm{Ai}_{-1}}$ in the
latter. Likewise for (\ref{eq61}) and (\ref{eq63}) below.


\section{L-G and Airy expansions: Case 1b}

Now let us consider the case $0<a_{0}\leq a^{2}\leq1-\delta$ ($-u\left(
1-{a_{0}}\right)  \leq\alpha\leq-u\delta<0$). As mentioned earlier, in this
case the expansions (\ref{eq51a}) and (\ref{eq43b}) remain valid for $U\left(
{n+\alpha+1,\alpha+1,uze^{-\pi i}}\right)  $. However, the same is not true
for the corresponding expansions (\ref{eq33}) and (\ref{eq43}) of
$L_{n}^{\left(  \alpha\right)  }\left(  {uz}\right)  $. This is because
$w_{0}\left(  {u,z}\right)  $ (defined by (\ref{eq31})) is in general not
recessive at $z=0$ when $\alpha<0$ (except when $\alpha$ is an integer).
Indeed, this time the appropriate solution recessive at the origin is given by%
\begin{equation}
w_{2}\left(  {u,z}\right)  \equiv z^{\left(  {\alpha+1}\right)  /2}%
e^{-uz/2}\mathbf{N}\left(  -n,\alpha+1,uz\right)  , \label{eq431}%
\end{equation}
where $\mathbf{N}\left(  a,c,x\right)  $ is Olver's scaled confluent
hypergeometric function defined in \cite[Chap. 7, sec. 9]{Olver:1997:ASF}. In
particular we have for ${\alpha\neq1,2,3,\cdots}$
\begin{equation}
\mathbf{N}\left(  -n,\alpha+1,x\right)  =\frac{x^{-\alpha}}{\Gamma\left(
{1-\alpha}\right)  }\left\{  {1+\mathcal{O}\left(  x\right)  }\right\}
\ \left(  {x\rightarrow0}\right)  , \label{eq432}%
\end{equation}
and hence%
\begin{equation}
w_{2}\left(  {u,z}\right)  \equiv\frac{z^{\left(  1-{\alpha}\right)  /2}%
}{\Gamma\left(  {1-\alpha}\right)  u^{\alpha}}\left\{  {1+\mathcal{O}\left(
z\right)  }\right\}  \ \left(  {z\rightarrow0}\right)  . \label{q433}%
\end{equation}
Comparing this to (\ref{eq32}) we see that $w_{2}\left(  {u,z}\right)  $ has
the desired recessive behaviour at $z=0$ when $\alpha<0$. We remark that when
$\alpha$ is a negative integer $w_{0}\left(  {u,z}\right)  $ and $w_{2}\left(
{u,z}\right)  $ are multiples of one another, since if $\alpha=-p$ for\ any
integer $p\in\left[  1,n\right)  $ we have from (\ref{eq0}) and (\ref{eq31})
\begin{equation}
w_{0}\left(  {u,z}\right)  =\frac{\left(  -u\right)  ^{p}z^{\left(
{p+1}\right)  /2}}{p!}\left\{  {1+\mathcal{O}\left(  z\right)  }\right\}
\ \left(  {z\rightarrow0}\right)  . \label{eq435}%
\end{equation}

By identifying $w_{2}\left(  {u,z}\right)  $ and $f^{-1/4}\left(  z\right)
V_{n,2}\left(  {u,\xi}\right)  $ in a similar manner to (\ref{eq32a}), we find
for $-u\left(  1-{a_{0}}\right)  \leq\alpha\leq-u\delta<0$ that
\begin{equation}%
\begin{array}
[c]{l}%
\mathbf{N}\left(  -n,\alpha+1,uz\right)  \sim\dfrac{1}{\Gamma\left(
1-{\alpha}\right)  }\left(  {\dfrac{\left\vert \alpha\right\vert }{u}}\right)
^{1/2}\left(  {\dfrac{u+\alpha}{u}}\right)  ^{u/2}\left\{  {\dfrac{\left(
{u+\alpha}\right)  e}{\alpha^{2}uz}}\right\}  ^{\alpha/2}\\
\\
\times\dfrac{1}{\left\{  {\left(  {z_{1}-z}\right)  \left(  {z_{2}-z}\right)
}\right\}  ^{1/4}}\exp\left\{  {\dfrac{1}{2}uz-u\xi+\sum\limits_{s=1}^{\infty
}{\left(  {-1}\right)  ^{s}\dfrac{\hat{{E}}_{s}\left(  z\right)  -\mu_{s}%
}{u^{s}}}}\right\}  ,
\end{array}
\label{eqNLG}%
\end{equation}
uniformly for $z\in D_{2}^{+}\cup D_{2}^{-}$.

Likewise, analogously to (\ref{eq43})%
\begin{equation}%
\begin{array}
[c]{l}%
\mathbf{N}\left(  -n,\alpha+1,uz\right)  \sim\dfrac{2\pi^{1/2}u^{1/6}}%
{\Gamma\left(  1-{\alpha}\right)  }\left(  {\dfrac{\left\vert \alpha
\right\vert }{u}}\right)  ^{1/2}\left(  {\dfrac{u+\alpha}{u}}\right)  ^{u/2}\\
\\
\times\left\{  {\dfrac{\left(  {u+\alpha}\right)  e}{\alpha^{2}uz}}\right\}
^{\alpha/2}\left\{  {\dfrac{\zeta}{\left(  {z_{1}-z}\right)  \left(  {z_{2}%
-z}\right)  }}\right\}  ^{1/4}\exp\left\{  {\dfrac{1}{2}uz+\sum\limits_{s=0}%
^{\infty}{\dfrac{\mu_{2s+1}}{u^{2s+1}}}}\right\} \\
\\
\times\left\{  {\mathrm{Ai}\left(  {u^{2/3}\zeta}\right)  A\left(
{u,z}\right)  +}\mathrm{Ai}^{\prime}{\left(  {u^{2/3}\zeta}\right)  B\left(
{u,z}\right)  }\right\}  ,
\end{array}
\label{eqNAiry}%
\end{equation}
uniformly for $z\in\mathbf{D}$.

To obtain the desired expansions for the Laguerre polynomials we use
(\ref{eqLMU}) and the connection formula%
\begin{equation}%
\begin{array}
[c]{l}%
\mathbf{M}\left(  -n,\alpha+1,z\right)  =\dfrac{e^{\alpha\pi i}n!}%
{\Gamma\left(  {n+\alpha+1}\right)  }\mathbf{N}\left(  -n,\alpha+1,z\right) \\
-\dfrac{\sin\left(  \pi\alpha\right)  n!}{\pi}e^{z}U\left(  {n+\alpha
+1,\alpha+1,ze^{-\pi i}}\right)  ,
\end{array}
\label{eq436}%
\end{equation}
to obtain%
\begin{equation}%
\begin{array}
[c]{l}%
L_{n}^{(\alpha)}\left(  {uz}\right)  =e^{\alpha\pi i}\mathbf{N}\left(
-n,\alpha+1,z\right) \\
\\
-\pi^{-1}\sin\left(  \pi\alpha\right)  \Gamma\left(  {n+\alpha+1}\right)
e^{z}U\left(  {n+\alpha+1,\alpha+1,ze^{-\pi i}}\right)  .
\end{array}
\label{eq437}%
\end{equation}
On inserting the expansions (\ref{eq51a}) and (\ref{eqNLG}) into (\ref{eq437})
yields the desired \mbox{L-G} expansion for $L_{n}^{\left(  \alpha\right)
}\left(  {uz}\right)  $ for $-u\left(  1-{a_{0}}\right)  \leq\alpha
\leq-u\delta<0$. Likewise, the corresponding Airy expansion comes from
(\ref{eq43b}), (\ref{eqNAiry}), and (\ref{eq437}).

\section{L-G and Airy expansions: Case 2}

We note that ${S\left(  z\right)  <0}$ for $z\in\left(  z_{2},\infty\right)
$, where $S\left(  z\right)  $ is given by (\ref{eq9}). Hence the same is true
for ${f^{1/2}\left(  z\right)  =S\left(  z\right)  /}\left(  2z\right)  ,$ and
with this in mind we now define the L-G variable by
\begin{equation}
\tilde{\xi}=-\int_{z_{2}}^{z}{f^{1/2}\left(  t\right)  dt,} \label{eq44}%
\end{equation}
so that $\tilde{\xi}\geq0$ for $z\geq z_{2}$. Explicit integration yields
\begin{equation}%
\begin{array}
[c]{l}%
\tilde{\xi}=-\dfrac{1}{2}\left(  {a^{2}+1}\right)  \ln\left\{  {z-S\left(
z\right)  -a^{2}-1}\right\}  -\dfrac{1}{2}S\left(  z\right)  +\max\left\{
a^{2},1\right\}  \ln\left(  {2a}\right) \\
\\
-\dfrac{1}{2}\left\vert {a^{2}-1}\right\vert \ln\left\{  {\dfrac{\left(
{a^{2}+1}\right)  z-\left(  {a^{2}-1}\right)  ^{2}+\left\vert {a^{2}%
-1}\right\vert S\left(  z\right)  }{z}}\right\}  .
\end{array}
\label{eq45}%
\end{equation}
We observe that $\tilde{\xi}{=}\xi-\min\left\{  a^{2},1\right\}  \pi i$, where
$\xi$ is the corresponding L-G variable from case 1.

The L-G asymptotic solution that we employ this time is given by the expansion%
\begin{equation}
\tilde{V}_{n,2}\left(  {u,\tilde{\xi}}\right)  =\exp\left\{  -{u\tilde{\xi
}+\sum\limits_{s=1}^{n-1}{\left(  -{1}\right)  ^{s}\frac{\hat{E}_{s}\left(
z\right)  }{u^{s}}}}\right\}  +\tilde{\varepsilon}_{n,2}\left(  {u,\tilde{\xi
}}\right)  , \label{eq52a}%
\end{equation}
where the coefficients $\hat{E}_{s}\left(  z\right)  $ are the same as in
Cases 1a and 1b.

In (\ref{eq52a}) $\tilde{\varepsilon}_{n,2}\left(  {u,\tilde{\xi}}\right)
=e^{-u{\tilde{\xi}}}\mathcal{O}\left(  {u^{-n}}\right)  $ uniformly in a
certain unbounded domain $\tilde{\Xi}_{2}$, with $e^{u{\tilde{\xi}}}%
\tilde{\varepsilon}_{n,2}\left(  {u,\tilde{\xi}}\right)  \rightarrow0$ as
$\operatorname{Re}{\tilde{\xi}}\rightarrow\infty$ in the domain. Let
$\tilde{\Delta}$ denote the domain in the $\tilde{\xi}$ plane corresponding to
$\Delta$ in the $\xi$ plane: thus it is $\Delta$ shifted by the factor
$-\min\left\{  a^{2},1\right\}  \pi i$. Then the domain $\tilde{\Xi}_{2}$
comprises the $\tilde{\xi}$ point subset of $\tilde{\Delta}$ for which there
is a path $\mathcal{\tilde{P}}_{2}$ (say) linking $\tilde{\xi}$ with
$\tilde{\alpha}_{2}=+\infty$ (corresponding to $z=+\infty+i0$) and having the
properties (i) $\mathcal{\tilde{P}}_{2}$ consists of a finite chain of $R_{2}$
arcs, and (ii) as $t$ passes along $\mathcal{\tilde{P}}_{2}$ from
$\tilde{\alpha}_{2}$ to $\xi$, $\operatorname{Re}\left(  ut\right)  $ is nonincreasing.

\begin{figure}[ptb]
\vspace*{0.8cm}
\centerline{\includegraphics[height=6cm,width=8cm]{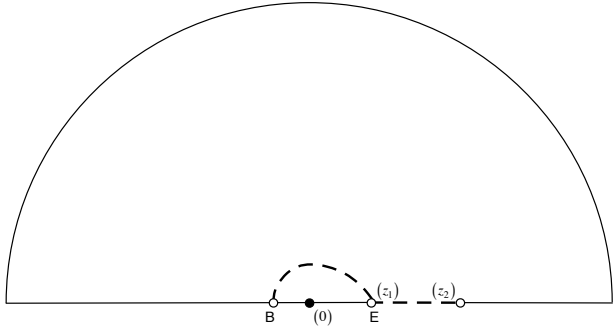}}\caption{Domain
$\tilde{D}_{2}^{+}$ in $z$ plane. The interval $\left[  z_{1},z_{2}\right]  $
and boundary EB are excluded.}%
\label{fig5zxi}%
\end{figure}

If $\tilde{\Xi}_{2}^{+}$ denotes the subset of $\tilde{\Xi}_{2}$ corresponding
to $0\leq\arg\left(  z\right)  \leq\pi$, then due to condition (ii) above
$\tilde{\Xi}_{2}^{+}$ must also exclude all points in the fourth quadrant of
the $\tilde{\xi}$ plane. The corresponding $z$-domain $\tilde{D}_{2}^{+}$
(say) is thus unshaded region depicted in Fig. \ref{fig5zxi}, where the
interval $\left[  z_{1},z_{2}\right]  $ and the boundary curve EB must be
excluded; the curve EB emanates from $z=z_{1}$ to a point on the upper part of
the cut along $\left(  -\infty,0\right)  $, and is given parametrically by
\begin{equation}
\int_{z_{1}}^{z}{f^{1/2}\left(  t\right)  dt=-\tau i,\ 0\leq\tau\leq\dfrac
{1}{2}\left\vert {a^{2}-1}\right\vert .}%
\end{equation}
We remark that $z=0$ is not contained in $\tilde{D}_{2}^{+}$.

We now use the behaviour of the Laguerre polynomials at infinity to match it
with the L-G expansion (\ref{eq52a}). Specifically, we know that
\begin{equation}
L_{n}^{\left(  \alpha\right)  }\left(  {uz}\right)  =\frac{\left(
{-1}\right)  ^{n}}{n!}U\left(  {-n,\alpha+1,uz}\right)  , \label{eq47}%
\end{equation}
and hence the following solution of (\ref{eq3})
\begin{equation}
\tilde{w}_{0}\left(  {u,z}\right)  \equiv\left(  {uz}\right)  ^{\left(
{\alpha+1}\right)  /2}e^{-uz/2}L_{n}^{\left(  \alpha\right)  }\left(
{uz}\right)  , \label{eq48}%
\end{equation}
has, from \cite[Eq. 13.7.3]{Daalhuis:2010:CHF}, the behaviour as
$z\rightarrow\infty$
\begin{equation}
\tilde{w}_{0}\left(  {u,z}\right)  =\frac{\left(  {-1}\right)  ^{n}}%
{n!}\left(  {uz}\right)  ^{\left(  {2u+\alpha}\right)  /2}e^{-uz/2}\left\{
{1+\mathcal{O}\left(  {\frac{1}{z}}\right)  }\right\}  . \label{eq49}%
\end{equation}
On comparing this to (\ref{eq51}) we see that $\tilde{w}_{0}\left(
{u,z}\right)  $ is the unique solution that is recessive in right half plane.
On matching, we therefore ascertain that there exits a constant $\tilde{c}%
_{0}\left(  u\right)  $ such that
\begin{equation}
\tilde{w}_{0}\left(  {u,z}\right)  =\tilde{c}_{0}\left(  u\right)
f^{-1/4}\left(  z\right)  \tilde{V}_{n,2}\left(  {u,}\tilde{\xi}\right)  ,
\label{eq52}%
\end{equation}
where $\tilde{V}_{n,2}\left(  {u,}\tilde{\xi}\right)  $ is given by
(\ref{eq52a}).

Using
\begin{equation}
\tilde{\xi}={\tfrac{1}{2}}z-{\tfrac{1}{2}}\left(  {a^{2}+1}\right)  \left\{
{\ln\left(  z\right)  +1}\right\}  +a^{2}\ln\left(  a\right)  +\mathcal{O}%
\left(  {z^{-1}}\right)  \ \left(  {z\rightarrow\infty}\right)  , \label{eq53}%
\end{equation}
we have from (\ref{eq52a}), (\ref{eq26}) and (\ref{eq49})
\begin{equation}%
\begin{array}
[c]{l}%
\tilde{c}_{0}\left(  u\right)  =\lim\limits_{z\rightarrow\infty}\left\{
{\dfrac{f^{1/4}\left(  z\right)  \tilde{w}_{0}\left(  {u,z}\right)  }%
{\tilde{V}_{n,2}\left(  {u,}\tilde{\xi}\right)  }}\right\} \\
\sim\dfrac{\left(  {-1}\right)  ^{n}u^{u/2}\left(  {u+\alpha}\right)
^{\left(  {u+\alpha}\right)  /2}}{\sqrt{2}e^{u+\left(  {\alpha/2}\right)  }%
n!}\exp\left\{  {-\sum\limits_{s=1}^{\infty}{\left(  {-1}\right)  ^{s}%
\dfrac{\lambda_{s}}{u^{s}}}}\right\}  .
\end{array}
\label{eq54}%
\end{equation}
Hence from (\ref{eq48}), (\ref{eq52}) and (\ref{eq54})
\begin{equation}%
\begin{array}
[c]{l}%
L_{n}^{\left(  \alpha\right)  }\left(  {uz}\right)  \sim\dfrac{\left(
-1\right)  ^{n}u^{\left(  u-1\right)  /2}\left(  u+\alpha\right)  ^{\left(
u+\alpha\right)  /2}}{e^{u+\left(  \alpha/2\right)  }n!\left(  uz\right)
^{\alpha/2}\left\{  {\left(  z-{z_{1}}\right)  \left(  z-{z_{2}}\right)
}\right\}  ^{1/4}}\\
\\
\times\exp\left\{  {\dfrac{1}{2}uz-u\tilde{\xi}+\sum\limits_{s=1}^{\infty
}{\left(  {-1}\right)  ^{s}\dfrac{\hat{{E}}_{s}\left(  z\right)  -\lambda_{s}%
}{u^{s}}}}\right\}  ,
\end{array}
\label{eq55}%
\end{equation}
uniformly for $z\in\tilde{D}_{2}^{+}\cup\tilde{D}_{2}^{-}$, where $\tilde
{D}_{2}^{-}$ is the conjugate region of $\tilde{D}_{2}^{+}$. We emphasise that
this expansion is not valid on the interval $\left[  z_{1},z_{2}\right]  $.

\begin{figure}[ptb]
\vspace*{0.8cm}
\centerline{\includegraphics[height=6cm,width=8cm]{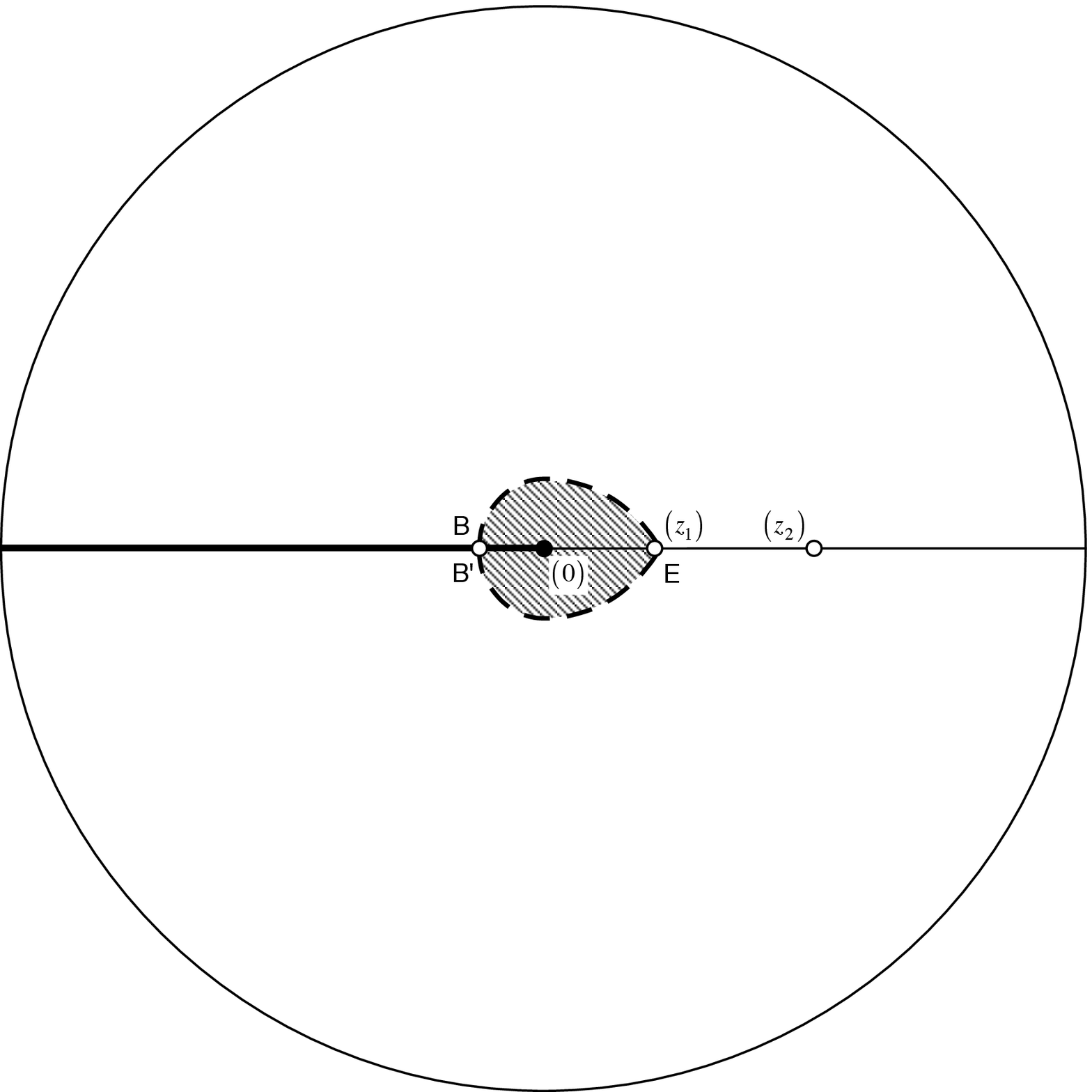}}\caption{Domain
$\mathbf{\tilde{D}}$ in $z$ plane.}%
\label{fig6zxi}%
\end{figure}

We next construct an Airy function expansion, which is valid at $z=z_{2}$,
similarly to (\ref{eq43}). We have, again by matching solutions (\ref{eq37})
with $j=0$ (all recessive at $z=+\infty$) with (\ref{eq48}), that
\begin{equation}%
\begin{array}
[c]{l}%
\tilde{C}_{0}\left(  u\right)  \left\{  {\dfrac{\left(  {z-z_{1}}\right)
\left(  {z-z_{2}}\right)  }{\tilde{\zeta}}}\right\}  ^{1/4}z^{\alpha
/2}e^{-uz/2}L_{n}^{\left(  \alpha\right)  }\left(  {uz}\right) \\
\\
=\mathrm{Ai}\left(  {u^{2/3}}\tilde{\zeta}\right)  \tilde{A}\left(
{u,z}\right)  +\mathrm{Ai}^{\prime}\left(  {u^{2/3}\tilde{\zeta}}\right)
\tilde{B}\left(  {u,z}\right)  ,
\end{array}
\label{eq56}%
\end{equation}
for some constant $\tilde{C}_{0}\left(  u\right)  $. Here $\tilde{\zeta
}=\left(  {{\tfrac{3}{2}}}\tilde{\xi}\right)  ^{2/3}$, where $\tilde{\xi}$ is
given by (\ref{eq45}), and the coefficient functions $\tilde{A}\left(
{u,z}\right)  \ $and $\tilde{B}\left(  {u,z}\right)  $ are analytic at
$z=z_{2}$ ($\tilde{\zeta}=0$). These have the expansions (c.f. (\ref{eq38})
and (\ref{eq39}))
\begin{equation}%
\begin{array}
[c]{l}%
\tilde{A}\left(  {u,z}\right)  \sim\exp\left\{  {\sum\limits_{s=1}^{\infty
}{\dfrac{\hat{{E}}_{2s}\left(  z\right)  +\tilde{{a}}_{2s}\tilde{\xi}%
^{-2s}/\left(  {2s}\right)  }{u^{2s}}}}\right\} \\
\times\cosh\left\{  {\sum\limits_{s=0}^{\infty}{\dfrac{\hat{{E}}_{2s+1}\left(
z\right)  -\tilde{{a}}_{2s+1}\tilde{\xi}^{-2s-1}/\left(  {2s+1}\right)
}{u^{2s+1}}}}\right\}  ,
\end{array}
\label{eq56a}%
\end{equation}
and
\begin{equation}%
\begin{array}
[c]{l}%
\tilde{B}\left(  {u,z}\right)  \sim\dfrac{1}{u^{1/3}\tilde{\zeta}^{1/2}}%
\exp\left\{  {\sum\limits_{s=1}^{\infty}{\dfrac{\hat{{E}}_{2s}\left(
z\right)  +a_{2s}\tilde{\xi}^{-2s}/\left(  {2s}\right)  }{u^{2s}}}}\right\} \\
\times\sinh\left\{  {\sum\limits_{s=0}^{\infty}{\dfrac{\hat{{E}}_{2s+1}\left(
z\right)  -a_{2s+1}\tilde{\xi}^{-2s-1}/\left(  {2s+1}\right)  }{u^{2s+1}}}%
}\right\}  .
\end{array}
\label{eq56b}%
\end{equation}
uniformly for $z\in\mathbf{\tilde{D}}$; this domain is the unshaded region
depicted in Fig. \ref{fig6zxi}, where EB' is the conjugate curve of EB. All
points on the boundary curve B'EB are excluded from $\mathbf{\tilde{D}}$.

Now as $\tilde{\zeta}\rightarrow\infty$ we find from (\ref{eq56a}),
(\ref{eq56b}), and the behaviour of the Airy function and its derivative for
large argument \cite[\S 9.7(ii)]{Olver:2010:AAR}, that the RHS behaves as
\begin{equation}%
\begin{array}
[c]{l}%
\mathrm{Ai}\left(  {u^{2/3}\tilde{\zeta}}\right)  \tilde{A}\left(
{u,z}\right)  +\mathrm{Ai}^{\prime}\left(  {u^{2/3}\tilde{\zeta}}\right)
\tilde{B}\left(  {u,z}\right) \\
\\
\sim\dfrac{1}{2\pi^{1/2}u^{1/6}\tilde{\zeta}^{1/4}}\exp\left\{  {-u\tilde{\xi
}-\sum\limits_{s=0}^{\infty}{\dfrac{\lambda_{2s+1}}{u^{2s+1}}}}\right\}  .
\end{array}
\label{eq57}%
\end{equation}
On the other hand, from (\ref{eq49}), the LHS has the behaviour
\begin{equation}%
\begin{array}
[c]{l}%
\tilde{C}_{0}\left(  u\right)  \left\{  {\dfrac{\left(  {z-z_{1}}\right)
\left(  {z-z_{2}}\right)  }{\tilde{\zeta}}}\right\}  ^{1/4}z^{\alpha
/2}e^{-uz/2}L_{n}^{\left(  \alpha\right)  }\left(  {uz}\right) \\
\\
\sim\dfrac{\left(  {-1}\right)  ^{n}\tilde{C}_{0}\left(  u\right)
u^{n}z^{\left(  {2n+\alpha+1}\right)  /2}e^{-uz/2}}{n!\tilde{\zeta}^{1/4}}.
\end{array}
\label{eq58}%
\end{equation}
Thus from (\ref{eq4}), (\ref{eq53}), (\ref{eq56}) - (\ref{eq58}) we can solve
for $\tilde{C}_{0}\left(  u\right)  $, and from this deduce that
\begin{equation}%
\begin{array}
[c]{l}%
L_{n}^{(\alpha)}\left(  uz\right)  \sim\dfrac{\left(  {-1}\right)  ^{n}%
2\pi^{1/2}u^{\left(  {u}/2\right)  -\left(  {1/3}\right)  }\left(  {u+\alpha
}\right)  ^{\left(  {u+\alpha}\right)  /2}}{n!\left(  uz\right)  ^{\alpha/2}%
}\\
\\
\times\left\{  {\dfrac{\tilde{\zeta}}{\left(  {z-z_{1}}\right)  \left(
{z-z_{2}}\right)  }}\right\}  ^{1/4}\exp\left\{  {\dfrac{1}{2}uz-u-\dfrac
{1}{2}\alpha+\sum\limits_{s=0}^{\infty}{\dfrac{\lambda_{2s+1}}{u^{2s+1}}}%
}\right\} \\
\\
\times\left\{  {\mathrm{Ai}\left(  {u^{2/3}\tilde{\zeta}}\right)  \tilde
{A}\left(  {u,z}\right)  +}\mathrm{Ai}^{\prime}{\left(  {u^{2/3}\tilde{\zeta}%
}\right)  \tilde{B}\left(  {u,z}\right)  }\right\}  .
\end{array}
\label{eq59}%
\end{equation}
\qquad

For the complementary solution we have, equivalent to (\ref{eq51a}), the
asymptotic expansion
\begin{equation}%
\begin{array}
[c]{l}%
U\left(  n+\alpha+1,\alpha+1,uze^{-\pi i}\right)  \sim\dfrac{\left(
{-1}\right)  ^{n+1}e^{\alpha\pi i+u+\left(  \alpha/2\right)  }}{u^{\left(
u+1\right)  /2}\left(  u+\alpha\right)  ^{\left(  u+\alpha\right)  /2}}\\
\times\dfrac{1}{\left(  uz\right)  ^{\alpha/2}\left\{  \left(  {z-z_{1}%
}\right)  \left(  {z-z_{2}}\right)  \right\}  ^{1/4}}\\
\times\exp\left\{  {-\dfrac{1}{2}uz+u\tilde{\xi}+\sum\limits_{s=1}^{\infty
}{\dfrac{\hat{{E}}_{s}\left(  z\right)  -\lambda_{s}}{u^{s}}}}\right\}  .
\end{array}
\label{eq61}%
\end{equation}

Similarly to the derivation of (\ref{eq59}), using
\begin{equation}%
\begin{array}
[c]{l}%
\mathrm{Ai}_{1}\left(  {u^{2/3}\tilde{\zeta}}\right)  \tilde{A}\left(
{u,z}\right)  +\mathrm{Ai}_{1}^{\prime}\left(  {u^{2/3}\tilde{\zeta}}\right)
\tilde{B}\left(  {u,z}\right) \\
\sim\dfrac{e^{\pi i/6}}{2\pi^{1/2}u^{1/6}\tilde{\zeta}^{1/4}}\exp\left\{
{u\tilde{\xi}+\sum\limits_{s=0}^{\infty}{\dfrac{\lambda_{2s+1}}{u^{2s+1}}}%
}\right\}  ,
\end{array}
\label{eq62}%
\end{equation}
we arrive at%
\begin{equation}%
\begin{array}
[c]{l}%
U\left(  {n+\alpha+1,\alpha+1,uze^{-\pi i}}\right)  \sim\dfrac{\left(
{-1}\right)  ^{n+1}2\pi^{1/2}e^{-\pi i/6}e^{\alpha\pi i}}{u^{1/3}\left\{
{u\left(  {u+\alpha}\right)  }\right\}  ^{\left(  {u+\alpha}\right)
/2}z^{\alpha/2}}\\
\\
\times\left\{  {\dfrac{\tilde{\zeta}}{\left(  {z-z_{1}}\right)  \left(
{z-z_{2}}\right)  }}\right\}  ^{1/4}\exp\left\{  {-\dfrac{1}{2}uz+u+\dfrac
{1}{2}\alpha-\sum\limits_{s=0}^{\infty}{\dfrac{\lambda_{2s+1}}{u^{2s+1}}}%
}\right\} \\
\\
\times\left\{  {\mathrm{Ai}_{1}\left(  {u^{2/3}\tilde{\zeta}}\right)
\tilde{A}\left(  {u,z}\right)  +}\mathrm{Ai}_{1}^{\prime}{\left(
{u^{2/3}\tilde{\zeta}}\right)  \tilde{B}\left(  {u,z}\right)  }\right\}  .
\end{array}
\label{eq63}%
\end{equation}

\section{Numerical results}

Here we illustrate the accuracy of the new expansions for Cases 1a and 2. We
concentrate on Laguerre polynomials with $n$ large and $\alpha$ non-negative,
but analogous results are available for negative $\alpha$, as well as 
for the complementary confluent hypergeometric functions.

\subsection{Case 1a}

In this case the expansions are valid in domains containing $z=0$ and
$z=z_{1}$. The relevant L-G expansion is given by (\ref{eq33}), and after
truncating after $N\geq1$ terms we have
\begin{equation}%
\begin{array}
[c]{l}%
L_{n}^{\left(  \alpha\right)  }\left(  {uz}\right)  =\dfrac{\Gamma\left(
{n+\alpha+1}\right)  }{n!\Gamma\left(  {\alpha+1}\right)  }\left(
{\dfrac{\alpha}{u}}\right)  ^{1/2}\left(  {\dfrac{u}{u+\alpha}}\right)
^{u/2}\\
\times\left\{  {\dfrac{\alpha^{2}}{\left(  {u+\alpha}\right)  uez}}\right\}
^{\alpha/2}\dfrac{1}{\left\{  {\left(  {z_{1}-z}\right)  \left(  {z_{2}%
-z}\right)  }\right\}  ^{1/4}}\\
\times\exp\left\{  {\dfrac{1}{2}uz-u\xi+\sum\limits_{s=1}^{N}{\left(
{-1}\right)  ^{s}\dfrac{\hat{{E}}_{s}\left(  z\right)  -\mu_{s}}{u^{s}}}%
}\right\}  \left\{  1+\mathcal{O}\left(  \dfrac{1}{u^{N+1}}\right)  \right\}
.
\end{array}
\label{eq63aa}%
\end{equation}
The order term is uniformly valid for $z\in D_{2}^{+}\cup D_{2}^{-}$, where
$D_{2}^{+}$ is shown in Fig. \ref{fig3zxi}, and $D_{2}^{-}$ is the conjugate
region of $D_{2}^{+}$.

The Airy expansion is given by (\ref{eq43}), and (\ref{eq38}) and (\ref{eq39})
are used to approximate the coefficient functions. Therefore uniformly for
$z\in\mathbf{D}$ (the domain of validity depicted in Fig. \ref{fig4zxi}) we
have%
\begin{equation}%
\begin{array}
[c]{l}%
L_{n}^{\left(  \alpha\right)  }\left(  {uz}\right)  =\dfrac{2\pi^{1/2}%
u^{1/6}\Gamma\left(  {n+\alpha+1}\right)  }{n!\Gamma\left(  {\alpha+1}\right)
}\left(  {\dfrac{\alpha}{u}}\right)  ^{1/2}\left(  {\dfrac{u}{u+\alpha}%
}\right)  ^{u/2}\\
\times\left(  {\dfrac{\alpha^{2}}{\left(  {u+\alpha}\right)  uez}}\right)
^{\alpha/2}\exp\left\{  {\dfrac{1}{2}uz+\sum\limits_{s=0}^{m-1}{\dfrac
{\mu_{2s+1}}{u^{2s+1}}}}\right\} \\
\times\left[  {\mathrm{Ai}\left(  {u^{2/3}\zeta}\right)  }\left\{
\mathcal{A}_{m}{\left(  {u,z}\right)  }+\mathcal{O}\left(  \dfrac{1}{u^{2m+1}%
}\right)  \right\}  \right. \\
{+}\left.  \mathrm{Ai}^{\prime}{\left(  {u^{2/3}\zeta}\right)  }\left\{
\mathcal{B}_{m}{\left(  {u,z}\right)  }+\mathcal{O}\left(  \dfrac
{1}{u^{2m+4/3}}\right)  \right\}  \right]  .
\end{array}
\label{eq63a}%
\end{equation}
In this, for any positive integer $m$, we have introduced the truncated
expansions%
\begin{equation}%
\begin{array}
[c]{l}%
\mathcal{A}_{m}\left(  {u,z}\right)  =\chi\left(  u,z\right)  \exp\left\{
{\sum\limits_{s=1}^{m}{\dfrac{\hat{{E}}_{2s}\left(  z\right)  +\tilde{{a}%
}_{2s}\xi^{-2s}/\left(  {2s}\right)  }{u^{2s}}}}\right\} \\
\times\cosh\left\{  {\sum\limits_{s=0}^{m-1}{\dfrac{\hat{{E}}_{2s+1}\left(
z\right)  -\tilde{{a}}_{2s+1}\xi^{-2s-1}/\left(  {2s+1}\right)  }{u^{2s+1}}}%
}\right\}  ,
\end{array}
\label{eq64}%
\end{equation}
and
\begin{equation}%
\begin{array}
[c]{l}%
\mathcal{B}_{m}\left(  {u,z}\right)  =\dfrac{\chi\left(  u,z\right)  }%
{u^{1/3}\zeta^{1/2}}\exp\left\{  {\sum\limits_{s=1}^{m}{\dfrac{\hat{{E}}%
_{2s}\left(  z\right)  +a_{2s}\xi^{-2s}/\left(  {2s}\right)  }{u^{2s}}}%
}\right\} \\
\times\sinh\left\{  {\sum\limits_{s=0}^{m-1}{\dfrac{\hat{{E}}_{2s+1}\left(
z\right)  -a_{2s+1}\xi^{-2s-1}/\left(  {2s+1}\right)  }{u^{2s+1}}}}\right\}  ,
\end{array}
\label{eq65}%
\end{equation}
where%
\begin{equation}
\chi\left(  u,z\right)  =\left\{  {\dfrac{\zeta}{\left(  {z_{1}-z}\right)
\left(  {z_{2}-z}\right)  }}\right\}  ^{1/4}. \label{eq65a}%
\end{equation}

We remark that the term $\chi\left(  u,z\right)  $ has been absorbed into the
scaled functions $\mathcal{A}_{m}\left(  {u,z}\right)  $ and $\mathcal{B}%
_{m}\left(  {u,z}\right)  $ because it has a removable singularity at
$z={z_{1}}$ ($\zeta=0$), and hence it is easier to compute via Cauchy's
integral formula when $z$ is close to ${z_{1}}$, as described next.

Let $z\in\mathbf{D}$: if this point is not too close to the turning point
$z_{1}$ we can use (\ref{eq64}) and (\ref{eq65}) directly for their
numerically stable computation in (\ref{eq63a}). On the other hand, if $z$ is
close to $z_{1}$ these expansions are not stable, since each $\hat{{E}}%
_{s}\left(  z\right)  $ is unbounded at this turning point. Instead we follow
\cite{Dunster:2017:COA} to compute these functions via Cauchy's integral
formula. Now, neither $\mathcal{A}_{m}{\left(  {u,z}\right)  }$ nor
$\mathcal{B}_{m}\left(  {u,z}\right)  $ are analytic at $z={z}_{1}$, but
$\chi\left(  u,z\right)  A\left(  u,z\right)  $ and $\chi\left(  u,z\right)
B\left(  u,z\right)  $ are. Thus, we have%
\begin{equation}%
\begin{array}{ll}
\chi\left(  u,z\right)  A\left(  u,z\right)  & = \dfrac{1}{2\pi i}
\displaystyle\oint_{\mathcal{L}_{1}}\dfrac{\chi\left(  u,t\right)  A\left(  u,t\right)
}{t-z}dt\\
& \\
&=\dfrac{1}{2\pi i}\displaystyle\oint_{\mathcal{L}_{1}}\dfrac{\mathcal{A}%
_{m}\left(  u,t\right)  +\mathcal{O}\left(  u^{-2m-1}\right)  }{t-z}dt\\
& \\
&=\dfrac{1}{2\pi i}\displaystyle\oint_{\mathcal{L}_{1}}\dfrac{\mathcal{A}_{m}\left(
u,t\right)  }{t-z}dt+\mathcal{O}\left(  \dfrac{1}{u^{2m+1}}\right)  ,
\end{array}
\label{eq67}%
\end{equation}
for some suitably-chosen bounded simple loop $\mathcal{L}_{1}$ in the
$t$-plane that is positively orientated, lies in an equivalent domain to
$\mathbf{D}$, and surrounds $t=z$ and $t=z_{1}$. Hence in (\ref{eq63a}) it is
legitimate (if necessary) to replace $\mathcal{A}_{m}{\left(  {u,z}\right)  }$
by $\left(  2\pi i\right)  ^{-1}\oint\limits_{\mathcal{L}_{1}}\mathcal{A}%
_{m}\left(  u,t\right)  \left(  t-z\right)  ^{-1}dt$, and similarly for
$\mathcal{B}_{m}\left(  {u,z}\right)  $. We then compute these integrals using
the trapezoidal rule, which has exponential convergence.

In our computations we shall take $N$ even and $m=N/2$ so that the order terms
in\ (\ref{eq63aa}) and (\ref{eq63a}) are of comparable magnitude (and likewise
for Case 2 considered later). Note that for fixed ${\zeta}$%
\begin{equation}
\begin{array}{l}
u^{1/6}{\mathrm{Ai}\left(  {u^{2/3}\zeta}\right)  =\exp}\left(  -{u\xi
}\right)  \mathcal{O}\left(  1\right)  ,\\
\\
\ u^{1/6}{\mathrm{Ai}}^{\prime
}{\left(  {u^{2/3}\zeta}\right)  =\exp}\left(  -{u\xi}\right)  \mathcal{O}%
\left(  u^{1/3}\right)  , \label{eq67a}%
\end{array}
\end{equation}
as $u\rightarrow\infty$ (see \cite[\S 9.7(ii)]{Olver:2010:AAR}).

If we consider the natural choice of $\mathcal{L}_{1}$ being a circular path
centred on $z_{1}$, it should have a radius smaller than
\begin{equation}
r_{m}=\min\left\{  z_{1},z_{2}-z_{1}\right\}  . \label{eq69}%
\end{equation}
We have that $r_{m}=z_{1}$ when $a\leq3+2\sqrt{2}$ (that is, $\alpha/\left(
n+1/2\right)  \leq4\left(  4+3\sqrt{2}\right)$) and $r_{m}=z_{2}-z_{1}$ if
$a\geq3+2\sqrt{2}$. When we have to use Cauchy's integral formula (for $z$
close to $z_{1}$) it is clear that $\operatorname{Re}z>0$; however, we have
shown that the validity of both the L-G approximation and the Airy expansion
away from the turning point $z_{1}$ extends to $\operatorname{Re}z<0$. We
first consider circular $\mathcal{L}_{1}$ as described above with a radius
smaller that $r_{m}$, and later we give details for $\operatorname{Re}z<0$
where we can use (\ref{eq64}) and (\ref{eq65}) directly (without resorting to
Cauchy integrals).

\begin{figure}[ptb]
\vspace*{0.8cm}
\centerline{\includegraphics[height=5.5cm,width=7.5cm]{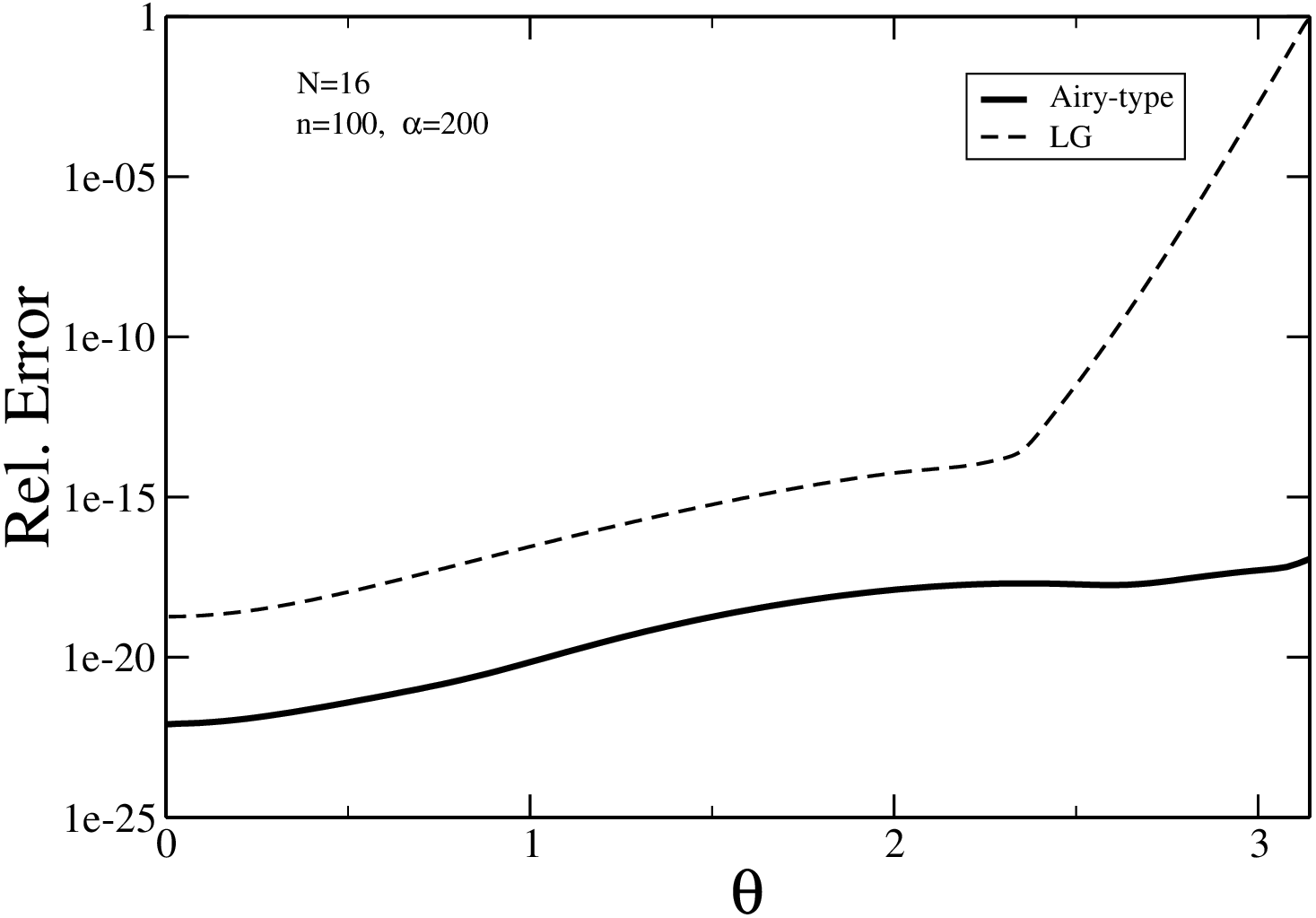}}\caption{Relative
errors $\epsilon_{rel}$ for the L-G expansion (\ref{eq63aa}) and the Airy
expansion (\ref{eq63a}), with coefficients computed from (\ref{eq64}) and
(\ref{eq65}), as a function of $\theta$, where $z=z_{1}-R e^{-i\theta}$,
$\theta\in[0,\pi]$ and $R=0.5 r_{m}=0.5 z_{1}$.}%
\label{fig5}%
\end{figure}

\begin{figure}[ptb]
\vspace*{0.7cm}
\centerline{\includegraphics[height=5.5cm,width=7.5cm]{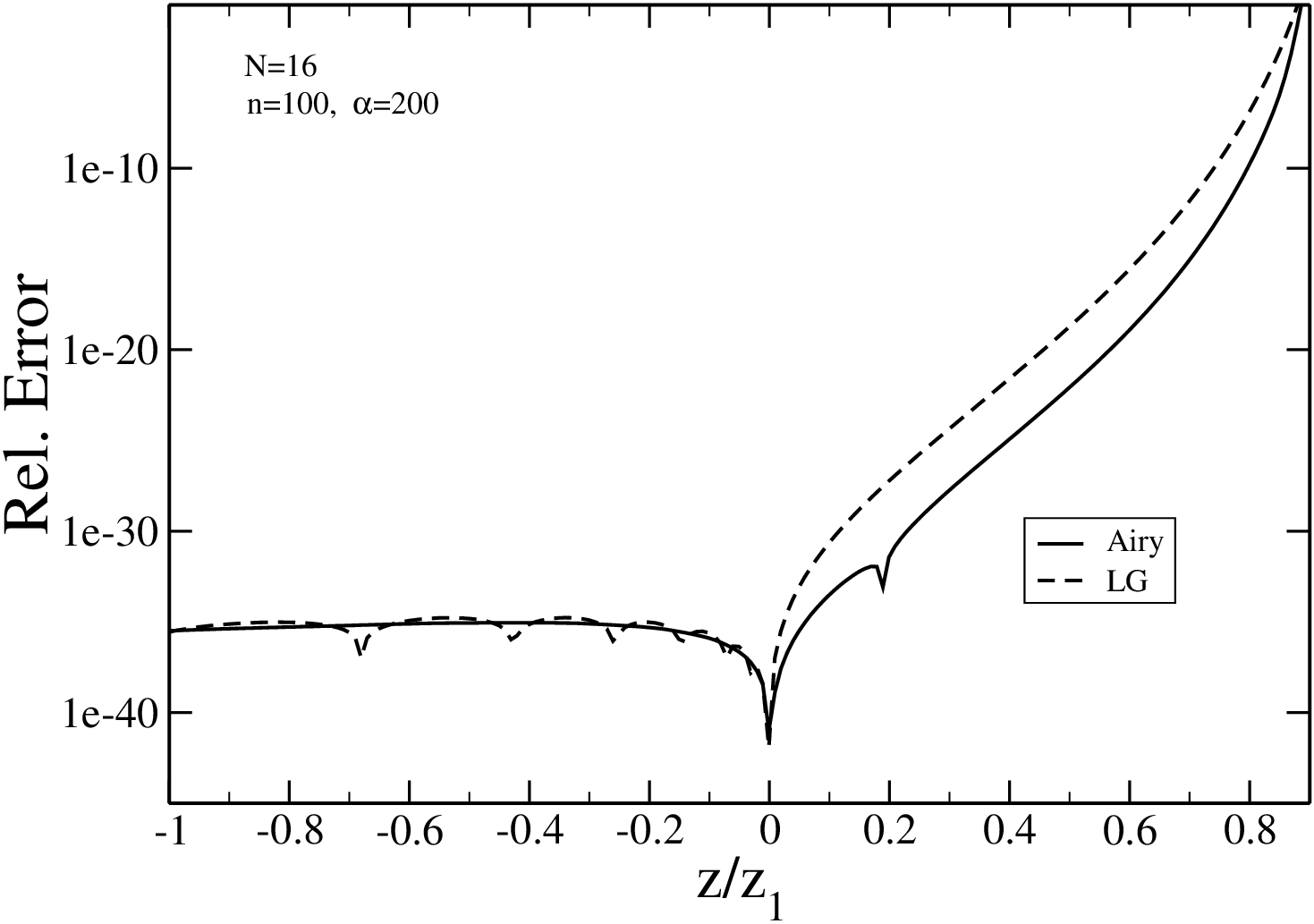}}\caption{Relative
errors $\epsilon_{rel}$ for the L-G expansion (\ref{eq63aa}) and the Airy
expansion (\ref{eq63a}), with coefficients computed from (\ref{eq64}) and
(\ref{eq65}), as a function of $z/z_{1}$ for real values of $z$, $z<z_{1}$.}%
\label{fig9}%
\end{figure}

Let the relative error of an approximation $f^{*}$ of a function $f$ be
defined in the usual manner by
\begin{equation}
\epsilon_{rel}=\left\vert \frac{f-f^{*}}{f}\right\vert . \label{eqRelErr}%
\end{equation}
Then in Fig. \ref{fig5} we show the relative errors of the L-G expansion
(\ref{eq63aa}) and the Airy expansion (\ref{eq63a}) for the Laguerre
polynomial. In both cases $\epsilon_{rel}$ is shown as a function of the angle
when the upper half circle centered at $z_{1}$ is followed clockwise. As
expected, we observe a step increase in the relative error for the L-G
approximation as we cross a Stokes line; contrarily, the Airy expansion works
well for all angles. In addition, we observe that the relative error for the
Airy expansion is smaller than for the L-G expansion for all $\theta$. Fig.
\ref{fig9} also shows that this is true for real values $0<z<z_{1}$, while for
negative $z$ both expansions give similar accuracies.

In these and in most of the figures, we consider the case of Laguerre
polynomials of degree $n=100$ and we use (\ref{eq63aa}) with $N=16$ and
(\ref{eq64}) and (\ref{eq65}) with $m=8$. If a higher degree $n$ or a larger
number of coefficients is considered, the accuracy will improve (except, of
course, if $N$ or $m$ is too large, because the series are asymptotic but not convergent).

\begin{figure}[ptb]
\vspace*{0.8cm}
\centerline{\includegraphics[height=5.5cm,width=7.5cm]{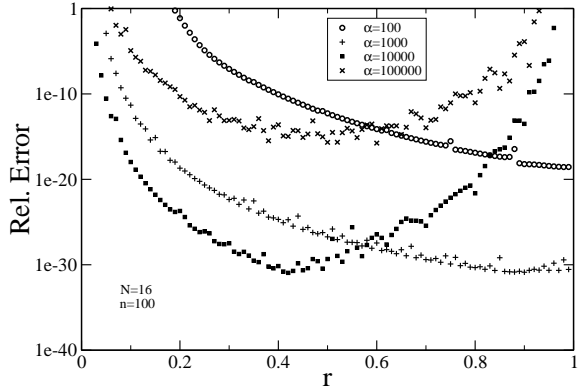}}\caption{The
maximum relative error of the Airy expansion (\ref{eq63a}), with coefficients
computed from (\ref{eq64}) and (\ref{eq65}), over the circles given by
$z=z_{2}+r r_{m} e^{i\theta}$, $\theta\in[0,2\pi)$, $r_{m}=\min\{z_{1}%
,z_{2}-z_{1}\}$, is plotted as a function of $r\in[0,1].$}%
\label{fig6}%
\end{figure}

Fig. \ref{fig6} shows the maximum relative error over circles centred at
$z_{1}$ of different radii (but always contained in the strip
$\{z:0<\operatorname{Re}z<z_{2}\}$). For the curves corresponding to
$a<3+2\sqrt{2}$ ($\alpha=100,1000$) the circles do not get close to $z_{2}$,
and this is why the relative error is monotonic as $r$ increases. In the other
to cases, the circle touches $z_{2}$ when $r=1$ which explains the different
pattern. We observe that the approximation works for $\alpha$ large comparable
to $n$ but that when $\alpha$ is too large accuracy is lost.

The circles considered in Fig. \ref{fig6} are the only circles centered at
$z_{1}$ that can be used for computing the coefficients of the Airy expansion
by Cauchy integrals. Indeed, we can not consider $z>z_{2}$ because the
expansions are meaningless there, and we can not consider $\operatorname{Re}%
z<0$ because of the discontinuity branch at the negative real axis (which
implies a discontinuity in the coefficients $A\left(  u,z\right)  $ and
$B\left(  u,z\right)  $).

\begin{figure}[ptb]
\vspace*{0.8cm}
\centerline{\includegraphics[height=5.5cm,width=7.5cm]{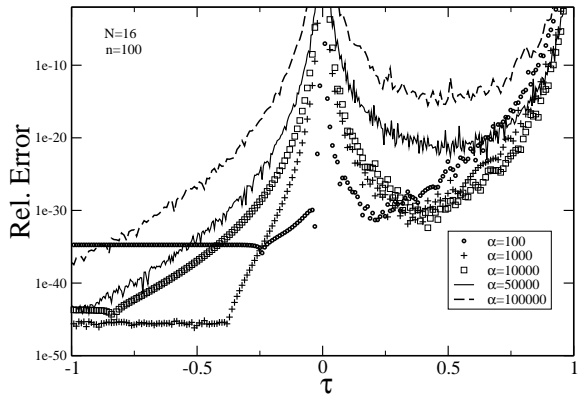}}\caption{Relative
error $\epsilon_{rel}$ of the Airy expansion (\ref{eq63a}), with coefficients computed from
(\ref{eq64}) and (\ref{eq65}), for real values of $z$ as a function of
$\tau=(z-z_{1})/(z_{2}-z_{1}).$}%
\label{fig7}%
\end{figure}

However, as commented earlier and shown in Fig. \ref{fig9}, the Airy expansion
(with coefficients computed by L-G asymptotics) is valid for
$\operatorname{Re}z<0$. In particular it is also valid for negative $z$, as
Fig. \ref{fig7} also shows. In this figure, we plot the relative error of the
Airy expansion (with coefficients from asymptotics) for real $z$ as a function
of $\tau=\left(  z-z_{1}\right)  /\left(  z_{2}-z_{1}\right)  $. This plot
includes negative values of $z$ for the two smaller values of $\alpha$, for
$\tau<\tau(0)\approx-0.0301$ when $\alpha=100$ and for $\tau<\tau
(0)\approx-0.4028$ when $\alpha=1000$ (for the other two cases $\tau(0)<-1$
and therefore $z>0$ in the figure). We observe that the expansion is also
accurate for negative $z$. As expected, the expansions fail close to the
turning points ($\tau=0$, $\tau=1$).

We notice that, as is well know, the Laguerre polynomials $L_{n}^{(\alpha
)}(uz)$ have $n$ positive real zeros when $\alpha>-1$ and most of them in the
interval $(z_{1},z_{2})$. Of course, the relative error (\ref{eqRelErr}) at
the zeros is meaningless, and loss of relative accuracy is unavoidable very
close to these zeros. In the previous figures, this loss of accuracy is not
clearly revealed because the function is sampled at values of $z$ which are
not too close to the zeros (the sample values are not selected to avoid the
zeros, simply happen to be not too close). Values of $z$ very close to the
zeros are needed to observe a significant accuracy loss; however, the previous
plots show some effect of the zeros because the relative error gives a
relatively noisy plot for values of $z$ where zeros occur when compared to the
cases without zeros. Compare, for example, the positive values of $\tau$ (for
which there are zeros) with the negative values (no zeros) in Fig. \ref{fig7}.
When we discuss Case 2 in the next section, we will show a detailed example of
computation close to the zeros.

\begin{figure}[ptb]
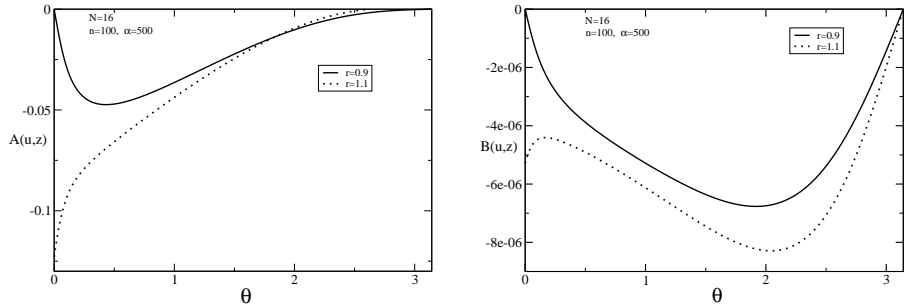

\vspace*{0.8cm}
\par
\begin{center}
\begin{minipage}{3cm}
\centerline{\includegraphics[height=5.6cm,width=5.6cm]{f8a.eps}}
\end{minipage}
\hspace*{3cm} \begin{minipage}{3cm}
\centerline{\includegraphics[height=5.6cm,width=5.6cm]{f8b.eps}}
\end{minipage}
\end{center}
\caption{Imaginary parts of the coefficient function approximations
$\mathcal{A}_{N/2}{\left(  {u,z}\right)  }$ and $\mathcal{B}_{N/2}{\left(
{u,z}\right)  }$ as a function of $\theta$, where $z=z_{1}-rz_{1}e^{-i\theta}$
for two values of $r$. }%
\label{fig8}%
\end{figure}

As commented before, even when the expansions are valid for $\operatorname{Re}%
z<0$, these values can not be used for computing the coefficients by means of
Cauchy integrals, due to the discontinuity at the branch cut. Fig. \ref{fig8}
shows this: we plot the imaginary part of the approximations $\mathcal{A}%
_{m}{\left(  {u,z}\right)  }$ and $\mathcal{B}_{m}{\left(  {u,z}\right)  }$ as
functions of the angle when the radius is such that the circle cuts the
negative real axis ($r=1.1$), and when it does not ($r=0.9$). For the second
case, both coefficients are real over the real line, but not in the first
case, when they have a non-zero imaginary part when $z$ is real and negative,
that is, when $\theta=0$ (we move along the circles clockwise).

Finally, we notice that the relative error close to the turning point becomes large, as 
shown in Fig. \ref{fig7}; this indicates that, as expected, the L-G expansions for the
coefficients of the Airy expansion are not accurate close to the turning
point. This loss of accuracy is repaired by computing these coefficients by Cauchy integrals
over a 
contour encircling the turning point (but sufficiently away from it) and
contained in the half-plane $\operatorname{Re}z>0$. Typically,
the error that can be obtained around the turning point is close to the
minimum error in Fig. \ref{fig6} (choosing as Cauchy contour a circle with the
radius corresponding to the minimum value of the error reached in that figure).

We give more explicit examples on the application of Cauchy integrals for the
Case 2 asymptotics we discuss next.

\subsection{Case 2}

Now we provide numerical evidence of the accuracy of the expansions for
Laguerre polynomials in domains containing the turning point $z_{2}$.

From (\ref{eq55}) the L-G approximation with $N$ terms is now given by
\begin{equation}%
\begin{array}
[c]{l}%
L_{n}^{\left(  \alpha\right)  }\left(  {uz}\right)  =\dfrac{\left(  -1\right)
^{n}u^{\left(  u-1\right)  /2}\left(  u+\alpha\right)  ^{\left(
u+\alpha\right)  /2}}{e^{u+\left(  \alpha/2\right)  }n!\left(  uz\right)
^{\alpha/2}\left\{  {\left(  z-{z_{1}}\right)  \left(  z-{z_{2}}\right)
}\right\}  ^{1/4}}\\
\times\exp\left\{  {\dfrac{1}{2}uz-u\tilde{\xi}+\sum\limits_{s=1}^{N}{\left(
{-1}\right)  ^{s}\dfrac{\hat{{E}}_{s}\left(  z\right)  -\lambda_{s}}{u^{s}}}%
}\right\}  \left\{  1+\mathcal{O}\left(  \dfrac{1}{u^{N+1}}\right)  \right\}
.
\end{array}
\label{eq70}%
\end{equation}
The order term is uniformly valid for $z\in\tilde{D}_{2}^{+}\cup\tilde{D}%
_{2}^{-}$, where $\tilde{D}_{2}^{+}$ is shown in Fig. \ref{fig5zxi}, and
$\tilde{D}_{2}^{-}$ is the conjugate region of $\tilde{D}_{2}^{+}$.

The Airy expansion is given by (\ref{eq43}), (\ref{eq56a}) and (\ref{eq56b}),
uniformly for $z\in\mathbf{\tilde{D}}$ (the domain of validity depicted in
Fig. \ref{fig6zxi}). Hence truncating similarly to (\ref{eq63a}) -
(\ref{eq65}) we have%
\begin{equation}%
\begin{array}
[c]{l}%
L_{n}^{\left(  \alpha\right)  }\left(  {uz}\right)  =\dfrac{\left(
{-1}\right)  ^{n}2\pi^{1/2}u^{\left(  {u}/2\right)  -\left(  {1/3}\right)
}\left(  {u+\alpha}\right)  ^{\left(  {u+\alpha}\right)  /2}}{n!\left(
uz\right)  ^{\alpha/2}}\\
\times\exp\left\{  {\dfrac{1}{2}uz-u-\dfrac{1}{2}\alpha+\sum\limits_{s=0}%
^{m-1}{\dfrac{\lambda_{2s+1}}{u^{2s+1}}}}\right\} \\
\times\left[  {\mathrm{Ai}\left(  {u^{2/3}\tilde{\zeta}}\right)  }\left\{
\widetilde{\mathcal{A}}_{m}{\left(  {u,z}\right)  }+\mathcal{O}\left(
\dfrac{1}{u^{2m+1}}\right)  \right\}  \right. \\
{+}\left.  \mathrm{Ai}^{\prime}{\left(  {u^{2/3}\tilde{\zeta}}\right)
}\left\{  \widetilde{\mathcal{B}}_{m}{\left(  {u,z}\right)  }+\mathcal{O}%
\left(  \dfrac{1}{u^{2m+4/3}}\right)  \right\}  \right]  ,
\end{array}
\label{eq71}%
\end{equation}
where, for any positive integer $m$,%
\begin{equation}%
\begin{array}
[c]{l}%
\widetilde{\mathcal{A}}_{m}\left(  {u,z}\right)  =\tilde{\chi}\left(
u,z\right)  \exp\left\{  {\sum\limits_{s=1}^{m}{\dfrac{\hat{{E}}_{2s}\left(
z\right)  +\tilde{{a}}_{2s}\tilde{\xi}^{-2s}/\left(  {2s}\right)  }{u^{2s}}}%
}\right\} \\
\times\cosh\left\{  {\sum\limits_{s=0}^{m-1}{\dfrac{\hat{{E}}_{2s+1}\left(
z\right)  -\tilde{{a}}_{2s+1}\tilde{\xi}^{-2s-1}/\left(  {2s+1}\right)
}{u^{2s+1}}}}\right\}  ,
\end{array}
\label{eq72}%
\end{equation}
and
\begin{equation}%
\begin{array}
[c]{l}%
\widetilde{\mathcal{B}}_{m}\left(  {u,z}\right)  =\dfrac{\tilde{\chi}\left(
u,z\right)  }{u^{1/3}\tilde{\zeta}^{1/2}}\exp\left\{  {\sum\limits_{s=1}%
^{m}{\dfrac{\hat{{E}}_{2s}\left(  z\right)  +a_{2s}\tilde{\xi}^{-2s}/\left(
{2s}\right)  }{u^{2s}}}}\right\} \\
\times\sinh\left\{  {\sum\limits_{s=0}^{m-1}{\dfrac{\hat{{E}}_{2s+1}\left(
z\right)  -a_{2s+1}\tilde{\xi}^{-2s-1}/\left(  {2s+1}\right)  }{u^{2s+1}}}%
}\right\}  ,
\end{array}
\label{eq73}%
\end{equation}
in which $\tilde{\chi}\left(  u,z\right)  =\left\{  {\dfrac{\tilde{\zeta}%
}{\left(  {z-z_{1}}\right)  \left(  {z-z_{2}}\right)  }}\right\}  ^{1/4}.$

As in Case 1a, if $z$ is close to the turning point (this time $z_{2}$) we can
replace $\widetilde{\mathcal{A}}_{m}{\left(  {u,z}\right)  }$ by $\left(  2\pi
i\right)  ^{-1}\oint\limits_{\mathcal{L}_{2}}\widetilde{\mathcal{A}}%
_{m}\left(  u,t\right)  \left(  t-z\right)  ^{-1}dt$ in (\ref{eq71}), and
similarly for $\widetilde{\mathcal{B}}_{m}\left(  {u,z}\right)  $. Here
$\mathcal{L}_{2}$ is a bounded contour in the $t$-plane that is positively
orientated, lies in an equivalent domain to $\mathbf{\tilde{D}}$, and
surrounds $t=z$ and $t=z_{2}$. As in Case 1a we typically take this to be a circle.

\begin{figure}[ptb]
\vspace*{0.8cm}
\centerline{\includegraphics[height=5.5cm,width=7.5cm]{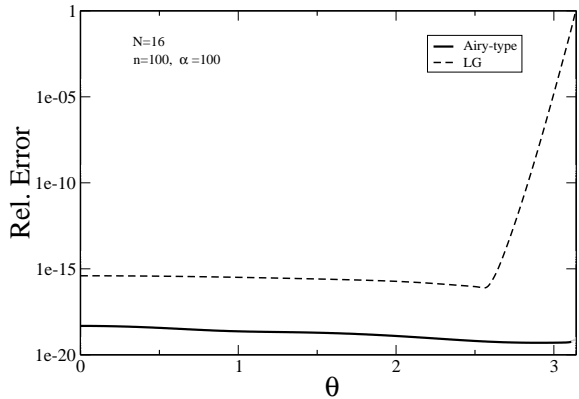}}\caption{Relative
errors $\epsilon_{rel}$ for the L-G expansion (\ref{eq70}) and the Airy
(\ref{eq71}) expansion, with coefficients computed from (\ref{eq72}) and
(\ref{eq73}), as a function of $\theta$, where $z=z_{2}+R e^{i\theta}$,
$\theta\in[0,\pi]$ and $R=0.25 (z_{2}-z_{1})$. }%
\label{fig1}%
\end{figure}

\begin{figure}[ptb]
\vspace*{0.8cm}
\centerline{\includegraphics[height=5.5cm,width=7.5cm]{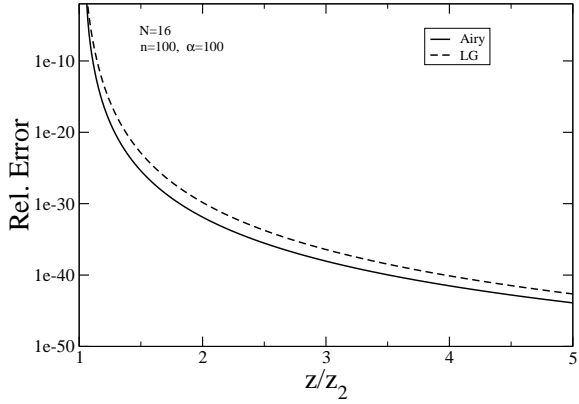}}\caption{Relative
errors $\epsilon_{rel}$ for the L-G expansion (\ref{eq70}) and the Airy
expansion (\ref{eq71}), with coefficients computed from (\ref{eq72}) and
(\ref{eq73}), as a function of $z/z_{2}$ for real values of $z$, $z>z_{2}$. }%
\label{fig10}%
\end{figure}

\begin{figure}[ptb]
\vspace*{0.8cm}
\centerline{\includegraphics[height=5.5cm,width=7.5cm]{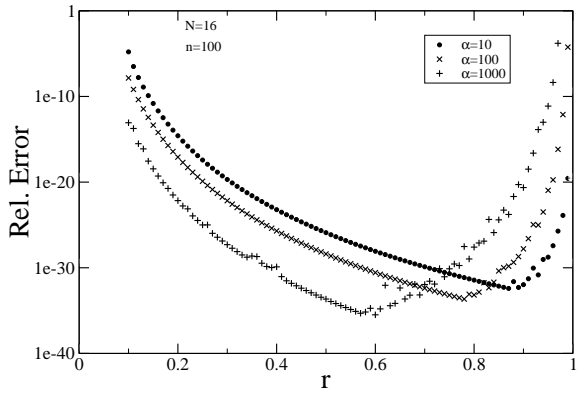}}\caption{The
maximum relative error of the Airy expansion (\ref{eq71}), with coefficients
computed from (\ref{eq72}) and (\ref{eq73}), over the circles given by
$z=z_{2}+Re^{i\theta}$, $\theta\in[0,2\pi)$ is plotted as a function of
$r=R/(z_{2}-z_{1})$. The circles are sampled with $100$ points. }%
\label{fig2}%
\end{figure}

\begin{figure}[ptb]
\vspace*{0.8cm}
\centerline{\includegraphics[height=5.5cm,width=7.5cm]{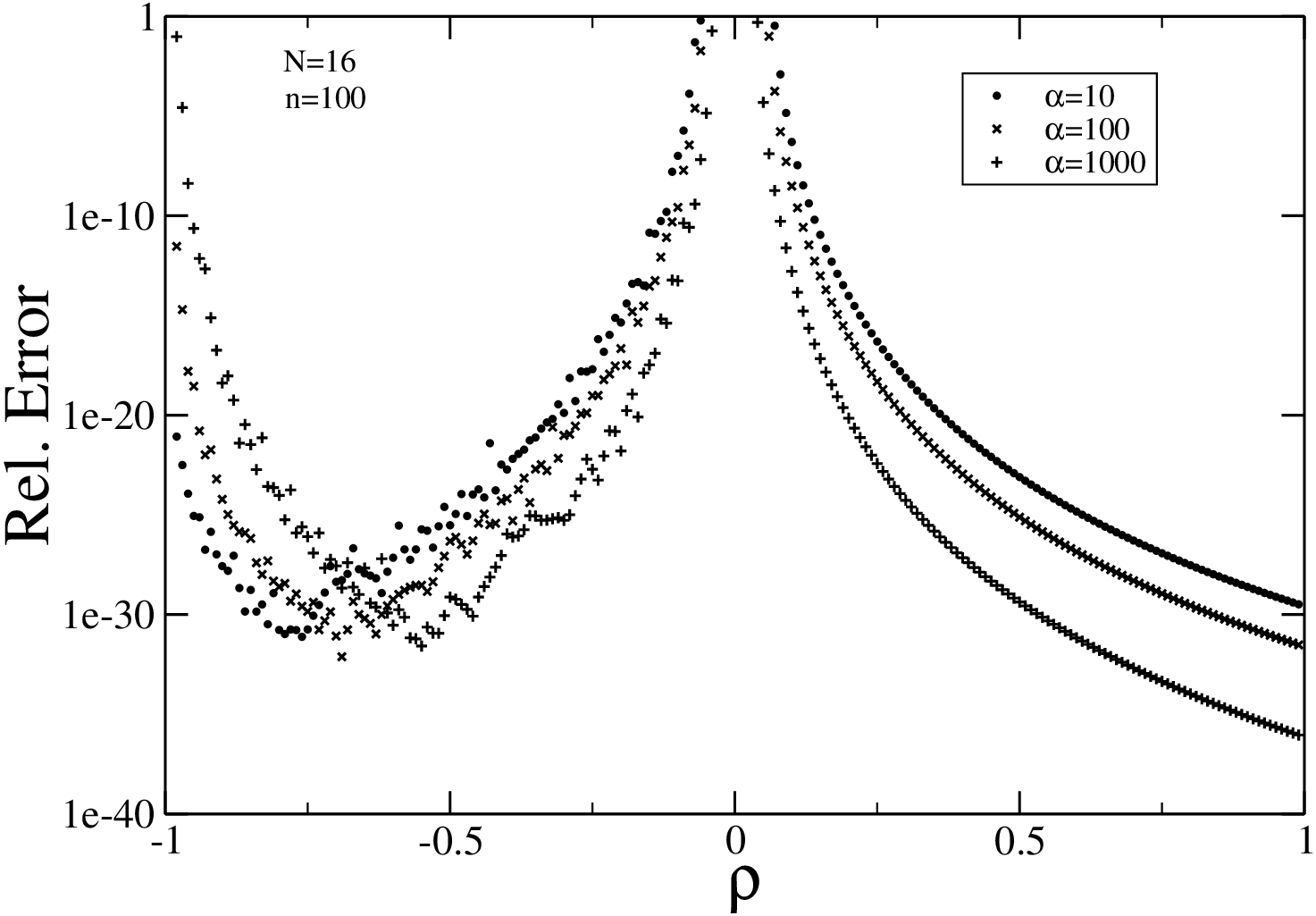}}\caption{Relative
error $\epsilon_{rel}$ of the Airy expansion (\ref{eq71}), with coefficients
computed from (\ref{eq72}) and (\ref{eq73}), for real values of $z$ as a
function of $\rho=\left(  z-z_{2}\right)  /\left(  z_{2}-z_{1}\right)  $. }%
\label{fig3}%
\end{figure}

\begin{figure}[ptb]
\vspace*{0.8cm}
\centerline{\includegraphics[height=5.5cm,width=7.5cm]{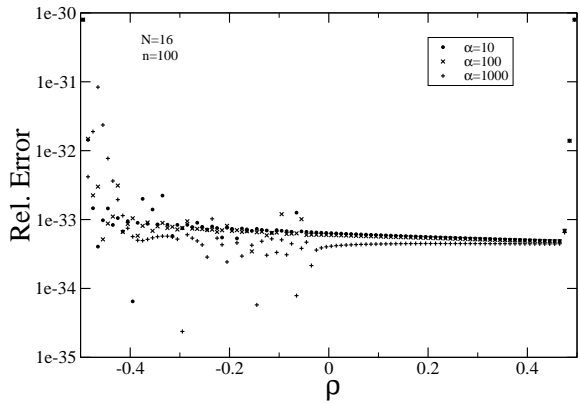}}\caption{Relative
error $\epsilon_{rel}$ of the Airy expansion (\ref{eq71}) for real values of
$z$ as a function of $\rho=\left(  z-z_{2}\right)  /\left(  z_{2}%
-z_{1}\right)  $, with the coefficient approximations $\widetilde{\mathcal{A}%
}_{m}{\left(  {u,z}\right)  }$ and $\widetilde{\mathcal{B}}_{m}{\left(
{u,z}\right)  }$ computed by Cauchy integrals. The Cauchy contour is a circle
of radius $R=0.7\left(  z_{2}-z_{1}\right)  $ centred at $z_{2}$.}%
\label{fig4}%
\end{figure}

Figs. \ref{fig1} and \ref{fig10} compare the accuracy attainable for the L-G
expansion (\ref{eq70}) and the Airy expansion (\ref{eq71}) - (\ref{eq73}),
both for real and complex variables. Figs. \ref{fig2} and \ref{fig3}
illustrate the accuracy of this same Airy expansion but for different values
of $\alpha$. And, finally, Fig. \ref{fig4} shows the accuracy of the Airy
expansion with the coefficients computed by Cauchy integrals. In all of these
cases, the coefficients of the Airy expansion are computed using $N=16$
($m=8$), and the degree of Laguerre polynomials is set to $n=100$.

In Fig. \ref{fig1}, we compare the L-G expansion with the Airy expansion away
from the turning point (with coefficient functions still computed with
(\ref{eq72}) and (\ref{eq73})). The relative accuracy over a circle of radius
$R=0.25\left(  z_{2}-z_{1}\right)  $ centered at $z_{2}$ is shown as a
function of the angle $\theta\in\lbrack0,\pi]$; this half-circle is circulated
counterclockwise. We observe that the L-G expansion (\ref{eq70}) tends to fail
for large $\theta$ in this interval, as the expansion loses meaning when we
cross a Stokes line. Contrarily, the Airy type expansion (\ref{eq71}) -
(\ref{eq73}) is valid in all the interval. In addition, we observe that the
relative error for the Airy expansion is smaller than for the L-G expansion
for all $\theta$. We further explore this fact for real variables $z>z_{2}$ in
Fig. \ref{fig10}.

Fig. \ref{fig2} shows the maximum relative error of the Airy expansion
(\ref{eq71}) - (\ref{eq73}), over circles centred at $z_{2}$. As expected, the
relative error increases both when the radius is small (because we are too
close to the turning point $z_{2}$) as well as when part of the circle becomes
too close to $z_{1}$.

Fig. \ref{fig3} shows again the relative error of the Airy expansion
(\ref{eq71}) - (\ref{eq73}) with coefficients computed with asymptotic series,
but for real variable. Again, the relative error increases close to the
turning points. For $z>z_{2}$ the relative error decreases as $z$ increases,
as can be expected. This figure shows that our alternative Cauchy integral
method for the computation of the coefficients $\widetilde{\mathcal{A}}%
_{m}{\left(  {u,z}\right)  }$ and $\widetilde{\mathcal{B}}_{m}\left(
{u,z}\right)  $ is needed around the turning point.

Fig. \ref{fig4} provides this computation using the trapezoidal rule with a
discretization of $100$ points in the upper half of the Cauchy contour (in the
lower half we consider complex conjugation). Combining this computation with
the results in Fig. \ref{fig3} we observe that it is possible to compute
accurately the Laguerre polynomials for $z>z_{1}$, but not too close to
$z_{1}$.

Observe that, according to Fig. \ref{fig2} (see also Fig. \ref{fig3}), over a
circle of radius $0.7(z_{2}-z_{1})$ centered at $z_{2}$, the relative errors
are of the order of $10^{-30}$, and we observe that the application of Cauchy
integrals permits us to maintain this accuracy inside the circle (but not too
close to the circle). As we see in Fig. \ref{fig4}, the relative error has
little variation and it is of the order of $10^{-30}$ inside the circle
$\left\vert z-z_{2}\right\vert <0.5\left(  z_{2}-z_{1}\right)  $ (the figure
is only for real $z$ but the same is true for complex $z$).

As we commented in the previous subsection, the \textit{relative} accuracy
unavoidably degrades very close to the zeros and this degration is, for
example, responsible for the different appearence of the graphs for $\rho>0$
and $\rho<0$ in Figs. \ref{fig3} and \ref{fig4}.

To illustrate the uniform accuracy of our approximations near the zeros we need to
replace the denominator of (\ref{eqRelErr}) with an "envelope" of the Laguerre
polynomial, which mimics the amplitude of $L_{n}^{\left(  \alpha\right)
}\left(  {x}\right)  $ but does not vanish at its zeros. For polynomials
having simple zeros, the envelope function $\mathrm{env}f(x) =\left\{
f^{2}(x)+f^{\prime2}(x)\right\} ^{1/2}$ serves this purpose. Then we define
the modified relative error of an approximation $f^{*}$ to a function $f$ as
\begin{equation}
\hat{\epsilon}_{rel}=\Frac{|f-f^*|}{\mathrm{env}f},
\end{equation}
where in our case $f(x)=L_{n}^{(\alpha)}(x)$.

\begin{figure}[ptb]
\vspace*{0.8cm}
\centerline{\includegraphics[height=5.5cm,width=7.5cm]{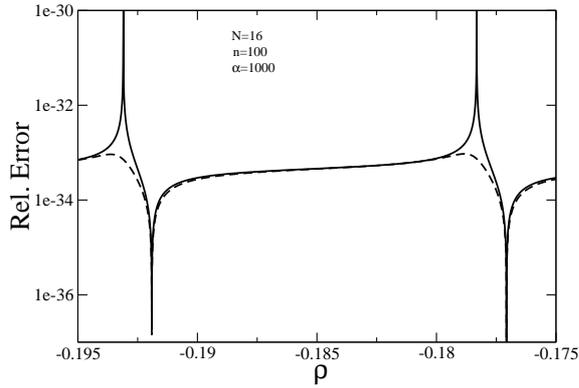}}\caption{Relative
error $\epsilon_{rel}$ (solid line) and modified relative error $\hat
{\epsilon}_{rel}$ (dashed line) of the Airy expansion (\ref{eq71}) for real
values of $z$ as a function of $\rho=\left(  z-z_{2}\right)  /\left(
z_{2}-z_{1}\right)  $, with the coefficient approximations $\widetilde
{\mathcal{A}}_{m}{\left(  {u,z}\right)  }$ and $\widetilde{\mathcal{B}}%
_{m}{\left(  {u,z}\right)  }$ computed by Cauchy integrals. The Cauchy contour
is a circle of radius $R=0.7\left(  z_{2}-z_{1}\right)  $ centered at $z_{2}%
$.}%
\label{errze}%
\end{figure}

In Fig. \ref{errze} we then show similar results for the Airy expansion in a
more restricted interval containing two zeros. In this both the relative error
$\epsilon_{rel}$ (solid line) and the modified relative error $\hat{\epsilon
}_{rel}$ (dashed line) are shown. In the figure we used many more sample
points than in Fig. \ref{fig4}, so that the (unavoidable) degradation of the
relative error $\epsilon_{rel}$ becomes more apparent, whereas in constrast
the modified relative error $\hat{\epsilon}_{rel}$ remains bounded. As can be
seen, the degradation of $\epsilon_{rel}$ only takes place very close to the
zeros. Of course, this \textit{relative} error degradation is unavoidable and
common to any method of numerical computation. On the other hand,
$\hat{\epsilon}_{rel}$ is small throughout the interval which illustrates the
uniform \textit{absolute} accuracy of our Airy expansion in the whole interval.

Finally, we give some additional results for other values of $n$, $\alpha$ and
$N$ as a further illustration of the accuracy of the computation using Cauchy
integrals. We pick a value of $z$ in the disc $\left\vert z-z_{2}\right\vert
<0.5\left(  z_{2}-z_{1}\right)  $. In particular, we fix $z=z_{1}%
+0.1(z_{2}-z_{1})$. We show the corresponding relative errors in Table
\ref{table1}. In the table, we take $300$ points in the upper half of the
Cauchy contour because with the previous selection ($100$ points) the
discretization error is not small enough in some cases; in particular for
$n=1000$ and $N=16$ ($m=N/2=8$), when the relative error becomes of the order
of $10^{-50}$.

\begin{table}[ptb]%
\[%
\begin{array}
[c]{|c|c|c|c|c|c|}\hline
&  &  &  &  & \\
N\Rightarrow & 2 & 4 & 8 & 12 & 16\\
\alpha,\, n &  &  &  &  & \\
\Downarrow &  &  &  &  & \\\hline
&  &  &  &  & \\
0,\,10 & 1.15\,10^{-6} & 2.88 \,10^{-9} & 2.20 \,10^{-13} & 1.28\,10^{-16} &
2.86\,10^{-19}\\\hline
&  &  &  &  & \\
0,\,100 & 1.46\,10^{-9} & 4.03 \,10^{-14} & 3.65 \,10^{-22} & 2.45\,10^{-29} &
6.24\,10^{-36}\\\hline
&  &  &  &  & \\
0,\,1000 & 1.50\,10^{-12} & 4.21 \,10^{-19} & 3.89 \,10^{-31} & 2.66\,10^{-42}
& 6.85\,10^{-53}\\\hline\hline
&  &  &  &  & \\
100,\,10 & 4.90\,10^{-7} & 1.22 \,10^{-9} & 9.53 \,10^{-13} & 5.68\,10^{-17} &
1.30\,10^{-19}\\\hline
&  &  &  &  & \\
100,\,100 & 7.79\,10^{-10} & 1.94 \,10^{-14} & 1.69 \,10^{-22} &
1.13\,10^{-29} & 2.88\,10^{-36}\\\hline
&  &  &  &  & \\
100,\,1000 & 1.31\,10^{-12} & 3.41 \,10^{-19} & 2.76 \,10^{-31} & 1.70
\,10^{-42} & 4.06 \,10^{-53}\\\hline
\end{array}
\]
\caption{Relative error $\epsilon_{rel}$ of the Airy expansion (\ref{eq71})
with $z=z_{2}+0.1(z_{2}-z_{1})$, where the coefficient approximations
$\widetilde{\mathcal{A}}_{m}{\left(  {u,z}\right)  }$ and $\widetilde
{\mathcal{B}}_{m}{\left(  {u,z}\right)  }$ are computed by Cauchy integrals
over the contour $|z-z_{2}|=0.7|z_{2}-z_{1}|$; different selections of the
degree $n$, the parameter $\alpha$ and the number of coefficients in the
asymptotic expansion $N=2m$ are considered.}%
\label{table1}%
\end{table}

We observe in Table \ref{table1} that the relative error shows little
variation with the value of $\alpha$, as the previous figures also showed, and
that the dependence on $n$ and $N$ is as expected for asymptotic
approximations with an error of $\mathcal{O}\left(  n^{-N-1}\right)  $.

This is shown more explicitly in Table \ref{table2}, where we give the values
of the computational asymptotic error constants $n^{N+1}\epsilon_{rel}$, with
$\epsilon_{rel}$ the relative errors in Table \ref{table1}; we observe these
experimental constants have a slow variation as a function of $n$ and $\alpha
$. While we also expect these constants to be $\mathcal{O}\left(  1\right)  $
they are in fact all quite small, which illustrates the uniform high accuracy
of our Airy expansions.

\begin{table}[ptb]%
\[%
\begin{array}
[c]{|c|c|c|c|c|c|}\hline
&  &  &  &  & \\
N\Rightarrow & 2 & 4 & 8 & 12 & 16\\
\alpha,\, n &  &  &  &  & \\
\Downarrow &  &  &  &  & \\\hline
&  &  &  &  & \\
0,\,10 & 1.15\,10^{-3} & 2.88 \,10^{-4} & 2.20 \,10^{-4} & 1.28\,10^{-3} &
2.86\,10^{-2}\\\hline
&  &  &  &  & \\
0,\,100 & 1.46\,10^{-3} & 4.03 \,10^{-4} & 3.65 \,10^{-4} & 2.45\,10^{-3} &
6.24\,10^{-2}\\\hline
&  &  &  &  & \\
0,\,1000 & 1.50\,10^{-3} & 4.21 \,10^{-4} & 3.89 \,10^{-4} & 2.66\,10^{-3} &
6.85\,10^{-2}\\\hline\hline
&  &  &  &  & \\
100,\,10 & 4.90\,10^{-4} & 1.22 \,10^{-4} & 9.53 \,10^{-5} & 5.68\,10^{-4} &
1.30\,10^{-2}\\\hline
&  &  &  &  & \\
100,\,100 & 7.79\,10^{-4} & 1.94 \,10^{-4} & 1.69 \,10^{-4} & 1.13\,10^{-3} &
2.88\,10^{-2}\\\hline
&  &  &  &  & \\
100,\,1000 & 1.31\,10^{-3} & 3.41 \,10^{-4} & 2.76 \,10^{-4} & 1.70 \,10^{-3}
& 4.06 \,10^{-2}\\\hline
\end{array}
\]
\caption{Computational asymptotic error constants estimated from the errors of
Table \ref{table1}.}%
\label{table2}%
\end{table}





\bibliographystyle{spmpsci}
\bibliography{biblio}

\end{document}